\documentclass[11pt]{amsart}
\newtheorem{theorem}{Theorem}[section]

\newtheorem{lemma}[theorem]{Lemma}

\theoremstyle{definition}
\newtheorem{definition}[theorem]{Definition}

\numberwithin{equation}{section}

\renewcommand{\tfrac}{\frac}
\newcommand{\subo}{0}
\newcommand{\frakH}{\mathfrak{H}}
\newcommand{\bofh}{\mathfrak{B}(\mathfrak{H})}

\newcommand{\empha}[1]{\textit{#1}}
\newcounter{numroman}

\begin{document}
\title{Local Unitary Cocycles of E$_0$-semigroups}
\author{Daniel Markiewicz}
\address{Daniel Markiewicz, Department of Mathematics,
Ben-Gurion University of the Negev, P.O.B. 653, Be'er Sheva 84105,
Israel.} \email{danielm@math.bgu.ac.il}
\author{Robert T. Powers}
\address{Robert T. Powers, Department of Mathematics,
University of Pennsylvania,
Phila\-delphia, PA 19104, U.S.A.}
\email{rpowers@math.upenn.edu}
\date{May 9, 2008}

\begin{abstract}
This paper concerns the structure of the group of local unitary cocycles, also called the gauge group, of an E$_0$-semigroup. The gauge group of a spatial E$_0$-semigroup has a natural action on the set of units by operator multiplication. Arveson has characterized completely the gauge group of E$_0$-semigroups of type I, and as a consequence it is known that in this case the gauge group action is transitive. In fact, if the semigroup has index $k$, then the gauge group action is transitive on the set of k+1-tuples of appropriately normalized independent units. An action of the gauge group having this property is called k+1-fold transitive.

We construct examples of E$_0$-semigroups of type II and index 1 which are not 2-fold transitive. These new examples also illustrate that an E$_0$-semigroup of type II$_k$ need not be a tensor product of an E$_0$-semigroup of type II$_0$ and another of type I$_k$.

\end{abstract}
\maketitle

\setcounter{section}{-1}
\section{Introduction}
%This paper is concerned with the structure of the gauge group of
%an E$_0$-semigroup of $\mathfrak{B}(\mathfrak{H})$, the algebra
%of all bounded operators on a separable Hilbert space $\frakH$.

An E$_0$-semigroup is a strongly continuous one-parameter semigroup of unit preserving $*$-endomorphisms of $\mathfrak{B}(\mathfrak{H})$, the algebra of
all bounded operators on a separable Hilbert space $\frakH$. In the 1930s Wigner showed that a one-parameter group of $*$-automorphisms of $\mathfrak{B}(\mathfrak{H})$  is always given by the action of a one-parameter strongly continuous unitary group by conjugation. In particular, the classification of one-parameter automorphism groups of $\bofh$ up to conjugacy can be reduced to the well-known multiplicity theory of Hahn-Helinger of unbounded self-adjoint operators. In contrast, the classification theory of E$_0$-semigroups up to cocycle conjugacy (which is the appropriate equivalence relation in this context) has proved to be much richer and full of surprises
(see for example \cite{Powers1987, Arveson1989a, Powers1999a, Tsirelson2000, Tsirelson2000a, Izumi2007a, Izumi2007, Izumi2007b}; we recommend Arveson's book \cite{A3} for an excellent introduction to the theory of E$_0$-semigroups).

One important cocycle conjugacy invariant of an E$_0$-semigroup is its gauge group (or more precisely the isomorphism class of the gauge group). Given an E$_0$-semigroup $\alpha$, a cocycle $C$ is a one-parameter strongly continuous family $\{C(t) : t \geq 0 \}$ satisfying the cocycle identity $C(t+s) = C(t) \alpha_t(C(s))$ for $t,s \geq 0$. The cocycle $C$ is said to be local if it satisfies the additional property that $C(t) \in
\alpha_t(\bofh)'$ for all $t\geq 0$. It is easy to verify that given two local cocycles $C_1$ and $C_2$, the expression $(C_1 \cdot C_2)(t) = C_1(t)C_2(t)$, for $t\geq 0$, defines another local cocycle. The set of unitary local cocycles is a group when endowed with this operation, and this group is called the gauge group of the E$_0$-semigroup $\alpha$. In terms of the product system approach to the study of E$_0$-semigroups, the gauge group is canonically isomorphic to the group of automorphisms of the product system associated with $\alpha$.

A unit for an E$_0$-semigroup $\alpha$ is a strongly continuous
one-parameter semigroup of isometries $\{U(t) : t \geq 0 \}$ which intertwines
the E$_0$-semigroup: $\alpha_t(A)U(t) = U(t) A$ for all $t\geq
0$, $A \in \bofh$. If an E$_0$-semigroups has at least one unit, it is called spatial. If it is spatial and it is generated by its units, it is called completely spatial or type I. All other spatial E$_0$-semigroups are called type II.  Non-spatial E$_0$-semigroups, also called type III, have been proven to exist by Powers \cite{Powers1987}, and in fact it follows from work of Tsirelson~\cite{Tsirelson2000a} that there exists a continuum of pairwise
non-cocycle conjugate examples. The type of an E$_0$-semigroup is also a cocycle conjugacy invariant.  We will only consider spatial semigroups in this paper, although we should note that, to our knowledge, little is known about the gauge group of a non-spatial E$_0$-semigroup.

Our main goal in this work is to study the action of the gauge group of a spatial E$_0$-semigroup on the set of units. If $U$ is a unit of $\alpha$ and $C$ is a local unitary cocycle, then $U'(t)=C(t)U(t)$ is another unit of $\alpha$, defining a natural action of the gauge group.

Completely spatial E$_0$-semigroups and their gauge groups are well-under\-stood. Every completely spatial E$_0$-semigroup is cocycle conjugate to a
CAR/CCR flow, and they are completely classified by the index
\cite{Arveson1989a, Powers1988, PowersPrice1990}. Furthermore, their gauge groups were completely characterized by
Arveson \cite{Arveson1989a,A3}.  One property which becomes apparent upon examining his characterization is that the gauge group of a completely spatial semigroup acts transitively on the set of units. In fact, even more can be gleaned from that characterization. If the completely spatial E$_0$-semigroup has index $k$, then any pair of $k+1$-tuples of appropriately normalized and independent units are related by an element of the gauge group. When the action of the gauge group of a spatial E$_0$-semigroup on its units has this
property, we say that the action is $k+1$-fold transitive.

Alevras, Price and Powers \cite{APP} were the first to break ground on the study of the gauge group of E$_0$-semigroups of type II. In their work, they characterize all contractive local cocycles (not just unitary local cocycles) for a certain class of  E$_0$-semigroups of type II and index zero. For semigroups of index zero, the set of units is essentially one-dimensional and the gauge group action is automatically ($1$-fold) transitive.

It is natural to inquire whether the action of the gauge group of a spatial E$_0$-semigroup of index $k$ on the units is always $k+1$-fold transitive.
In this paper we show that this is not the case, by constructing a class of E$_0$-semigroups of type II and index 1 whose gauge group action on the set of units is not $2$-fold transitive.

It is possible, although we could not verify it, that within the class which we have constructed there could also be examples of E$_0$-semigroups whose gauge group action on the set of units is not transitive. While this article was in preparation, it came to our attention that non-transitive examples were obtained by Tsirelson \cite{Tsirelson2004}, using different techniques. We do not know the exact relationship between his examples and our own. Nevertheless,  in the last section we discuss some features which they have in common.

We also observe that our examples, as well as Tsirelson's \cite{Tsirelson2004},
provide a direct answer to an old question. When Arveson \cite{Arveson1989}
proved that the index is additive with respect to tensor products, it was natural
to inquire whether type II$_k$ semigroups can be decomposed as
tensor products of type II$_0$ and I$_k$. Alas, that is not the
case, as the E$_0$-semigroups which we construct are of type
II$_1$ yet they are not tensor products of type II$_0$ and type
I$_1$ semigroups.

Our approach involves a detailed analysis of the E$_0$-semigroups obtained via minimal dilation of certain  CP-flows. Bhat \cite{B1} proved that CP-semigroups can be dilated to E$_0$-semigroups, and this result proved very useful for the construction and analysis of new examples of E$_0$-semigroups  (a very incomplete list of work in this direction includes \cite{Arveson1999, BhatSkeide2000, Markiewicz2003, MuhlySolel2002, Powers1999b}).  Bhat~\cite{Bh} has also found a one-to-one correspondence between the compressions of a CP-semigroup and the compressions of its minimal dilation. Pursuing this correspondence, Powers~\cite{BP} subsequently carried out a study of a class of CP-semigroups, called CP-flows, and their minimal dilations to E$_0$-semigroups. In particular, several results were obtained in \cite{BP} for the analysis of the cocycle conjugacy of the minimal dilations of CP-flows, as well as their contractive local cocycles. We make full use of this favorable framework, which is in fact quite general, given that all spatial E$_0$-semigroups arise from the minimal dilation of an appropriate CP-flow (see \cite{BP}).

We now provide an outline of the contents of the following sections. In Section 2 we describe in detail the basic background and terminology, with an emphasis on the material related to \cite{BP}. In Section 3 we introduce the class of examples which will be of interest, and describe some of its key properties. In Section 4 we turn to the analysis of the local cocycles of the E$_0$-semigroups under consideration. Finally, in the last section we summarize our main results.

\subsection*{Acknowledgments}
This work was partially carried out while D.M. was a Lecturer at the University of Pennsylvania and later a Technion Swiss Society Postdoctoral Fellow at the Technion in Haifa, Israel. He would like to thank R.T.P. and his wife for their wonderful generosity and hospitality during his time in Philadelphia. He also thanks Baruch Solel at the Technion, for his hospitality and many interesting conversations, and the friendly staff of both institutions.

\section{Background, Notation and Definitions}\label{section:background}

We begin with the definition of $E_0$-semigroups of $\mathfrak{B}(\mathfrak{H})$
the set of all bounded operators on a separable Hilbert space $\mathfrak{H}$.  For a
detailed discussion of $E_0$-semigroups we refer to Arveson's excellent
book \cite{A3}.

\begin{definition} We say $\alpha$ is an \empha{E$_0$-semigroup} of
$\mathfrak{B} (\mathfrak{H} )$ if the following conditions are satisfied.

\begin{enumerate}
%\widestnumber\item{(iii)} \roster

\item $\alpha_t$ is a $*$-endomorphism of $\mathfrak{B} (\mathfrak{H} )$
for
each $t \geq 0$.

\item $\alpha_0$ is the identity endomorphism and
$\alpha_t\circ\alpha _s = \alpha_{t+s}$ for all $s,t \geq 0$.

\item For each $\rho \in \mathfrak{B} (\mathfrak{H} )_*$ (the predual of
$\mathfrak{B} (\mathfrak{H} ))$ and $A \in \bofh$ the function $\rho (\alpha_t(A))$
is a continuous function of $t$.

\item $\alpha_t(I) = I$ for each $t \geq 0$ ($\alpha_t$ preserves the
unit).
\end{enumerate}
\end{definition}

The appropriate notions of when two $E_0$-semigroups are similar are
conjugacy and cocycle conjugacy (which comes from Alain Connes' definition
of outer conjugacy).

\begin{definition}  Suppose $\alpha$ and $\beta$ are
$E_0$-semigroups $\mathfrak{B} (\mathfrak{H}_1)$ and $\mathfrak{B}
(\mathfrak{H}_2)$. We say
$\alpha$ and $\beta$ are \empha{conjugate}, denoted $\alpha \approx \beta$, if
there is a $*$-isomorphism $\phi$ of $\mathfrak{B} (\mathfrak{H}_1)$ onto
$\mathfrak{B}
(\mathfrak{H}_2)$ so that $\phi\circ \alpha_t = \beta_t\circ\phi$ for all $t
\geq 0$.  We say $\alpha$ and $\beta$ are \empha{cocycle conjugate}, denoted
$\alpha_t \sim \beta_t$, if $\alpha ^{\prime}$ and $\beta$ are conjugate
where $\alpha$ and $\alpha ^{\prime}$ differ by a unitary cocycle (i.e.,
there is a strongly continuous one parameter family of unitaries $U(t)$ on
$\mathfrak{B} (\mathfrak{H}_1)$ for $t \geq 0$ satisfying the cocycle condition
$U(t)\alpha_t(U(s)) = U(t+s)$ for all $t,s \geq$ 0 so that
$\alpha_t^{\prime}(A) = U(t)\alpha_t(A)U(t)^{-1}$ for all $A \in \mathfrak{B}
(\mathfrak{H}_1)$ and $t \geq 0)$.
\end{definition}

An $E_0$-semigroup $\alpha_t$ is \empha{spatial} if there is a semigroup of
isometries $U(t)$ which intertwine so $U(t)A = \alpha_t(A)U(t)$ for $A \in
\mathfrak{B} (\mathfrak{H} )$ and $t > 0$.  The property of being spatial is a
cocycle conjugacy invariant.

An extremely useful and well known result in the theory of $C^*$-algebras
is the Gelfand-Segal construction of a cyclic $*$-representation of a
$C^*$-algebra associated with a state of the $C^*$-algebra.  In the study
of $E_0$-semigroups there is a result in the same spirit which says that
every  semigroup of unital completely positive maps of $\mathfrak{B}
(\mathfrak{K} )$
can be dilated to an $E_0$-semigroup of $\mathfrak{B} (\mathfrak{H} )$ where
$\mathfrak{H}$ can be thought of as a larger Hilbert space containing
$\mathfrak{K}$. We begin with a review of the properties of completely positive
maps.

A linear map $\phi$ from a $C^*$-algebra $\mathfrak{A}$ into $\mathfrak{B}
(\mathfrak{H})$ is \empha{completely positive} if
$$
\sum_{i,j=1}^n (f_i,\phi (A_i^*A_j)f_j) \geq 0
$$
for $A_i \in \mathfrak{A} ,\medspace f_i \in \mathfrak{H}$ for $i = 1,2,\dots
,n$ and
$n = 1,2,\dots $. Stinespring's central result is that if $\mathfrak{A}$ has a
unit and $\phi$ is a completely positive map from $\mathfrak{A}$ into
$\mathfrak{B}
(\mathfrak{H} )$ then there is a $*$-representation $\pi$ of $\mathfrak{A}$ on
$\mathfrak{B} (\mathfrak{K} )$ and an operator $V$ from $\mathfrak{H}$ to
$\mathfrak{K}$ so that $\phi
(A) = V^*\pi (A)V$ for $A \in \mathfrak{A} $.  And $\pi$ is determined by $\phi$
up to unitary equivalence if the linear span of $\{\pi (A)Vf\}$ for $A \in
\mathfrak{A}$ and $f \in \mathfrak{H}$ is dense in $\mathfrak{K} $.

Often we speak of one functional or map dominating another.  We introduce a
word for the functional or map that is dominated.  The word is
\empha{``subordinate"}.  If $A$ is an object which is positive with respect to some
order structure we say $B$ is a subordinate of $A$ if $B$ is the same kind
of thing $A$ is and $B$ is positive and $B$ is less than $A$.  For example
if we are speaking of the positive integers the subordinates of 4 are
4,3,2,1.  If $A$ is a positive operator then the subordinates of A are
operators $B$ with $A \geq B \geq 0$.  Suppose $E$ is a projection.  Are
the subordinates of a projection $E$ projections under $E$ or the operators
under $E$?  The answer depends on the context.

A \empha{$CP$-semigroup} of $\mathfrak{B} (\mathfrak{H} )$ is a strongly continuous one parameter semigroup of completely positive maps of $\mathfrak{B} (\mathfrak{H})$
into itself.  We now state Bhat's theorem \cite{B1} for $\mathfrak{B}
(\mathfrak{H}
)$.

\begin{theorem}\label{thm:Bhat's-dilation-thm} %{Theorem 1.3}
Suppose $\alpha$ is a unital
$CP$-semigroup of
$\mathfrak{B} (\mathfrak{H} )$. Then there is an $E_0$-semigroup $\alpha^d$ of
$\mathfrak{B}
(\mathfrak{K} )$ and an isometry $W$ from $\mathfrak{H}$ to $\mathfrak{K}$ so
that
$$
\alpha_t(A) = W^*\alpha_t^d(WAW^*)W
$$
and $\alpha_t(WW^*) \geq WW^*$ for $t > 0$ and if the projection $E = WW^*$
is minimal, which means the span of the vectors
$$
\alpha_{t_1}^d(EA_1E)\alpha_{t_2}^d(EA_2E)\cdots\alpha_{t_n}^d(EA_nE)Wf
$$
for $f \in \mathfrak{H} ,\medspace A_i \in \mathfrak{B} (\mathfrak{H}
),\medspace t_i \geq
0$ for $i = 1,2,\dots $ and $n = 1,2,\dots $ is dense in $\mathfrak{H}$, then
$\alpha^d$ is determined up to conjugacy.
\end{theorem}

We use Arveson's definition of minimality which is easier to state and
equivalent to Bhat's.

Suppose $\alpha$ is an $E_0$-semigroup of $\mathfrak{B} (\mathfrak{H} )$.  We
characterize the subordinates of $\alpha $, (i.e. the $CP$-semigroups
$\beta$ of $\mathfrak{B} (\mathfrak{H} )$ so that the mapping $A \rightarrow
\alpha_t(A) - \beta_t(A)$ is completely positive for all $t \geq 0)$. The
subordinates of $\alpha$ are given by positive local cocycles.  A cocycle
is a $\sigma$-weakly continuous one parameter family of operators $C(t)$
satisfying the cocycle relation
$$
C(t + s) = C(t)\alpha_t(C(s))
$$
for all $s,t \geq 0$.  The cocycle $C(t)$ is local if $C(t) \in
\alpha_t(\mathfrak{B} (\mathfrak{H} ))^{\prime}$ for all $t > 0$.  The local
cocycles
and their order structure are a cocycle conjugacy invariant.  As first
shown by Bhat \cite{B1} there is an order isomorphism from the subordinates
of a unital $CP$-semigroup of $\mathfrak{B} (\mathfrak{H} )$ to the subordinates
of
its minimal dilation to an $E_0$-semigroup of $\mathfrak{B} (\mathfrak{K} )$ We
use
the notation of \cite{BP}.

\begin{theorem}\label{thm:arveson's-minimal-dilation} %Theorem 1.4
Suppose $\alpha$ is a unital $CP$-semigroup of
$\mathfrak{B} (\mathfrak{H} )$ and $\alpha^d$ is the minimal dilation of
$\alpha$ to
an $E_0$-semigroup of $\mathfrak{B} (\mathfrak{K} )$ and $W$ is an isometry from
$\mathfrak{H}$ to $\mathfrak{K}$ so that $WW^*$ is a minimal projection for
$\alpha^d$ and
$$
\alpha_t(A) = W^*\alpha_t^d(WAW^*)W
$$
for $A \in \mathfrak{B} (\mathfrak{H} )$ and $t \geq 0$.  Then there is an order
isomorphism from the subordinates of $\alpha$ to the subordinates of
$\alpha^d$ given as follows.  Suppose $\gamma$ is a subordinate of
$\alpha^d$ and $C(t) = \gamma_t(I)$ for $t \geq 0$ is the local cocycle
associated with $\gamma$ then the subordinate $\beta$ of $\alpha$ under
this isomorphism is given by
$$
\beta_t(A) = W^*C(t)\alpha_t^d(WAW^*)W
$$
for $A \in \mathfrak{B} (\mathfrak{H} )$ and $t \geq 0$.
\end{theorem}

In this paper we will frequently make use of corners.  This is a trick
introduced by A. Connes.

\begin{definition}\label{def:corner}   %Definition 1.5
Suppose $\alpha$ and $\beta$ are
$CP$-semigroups of $\mathfrak{B} (\mathfrak{H} )$ and $\mathfrak{B}
(\mathfrak{K} )$. Then
$\gamma$ is a \empha{corner from $\alpha$ to $\beta$} if $\Theta$ given by
$$
\Theta_t \left[\begin{matrix} A&B
\\
C&D
\end{matrix} \right] =
\left[\begin{matrix} \alpha_t(A)&\gamma_t(B)
\\
\gamma_t^*(C)&\beta_t(D)
\end{matrix} \right]
$$
for $t \geq 0$ and $A \in \mathfrak{B} (\mathfrak{H} ),\medspace D \in
\mathfrak{B} (\mathfrak{K} ),\medspace B$ a linear operator from $\mathfrak{K}$
to $\mathfrak{H}$ and $C$ a
linear operator from $\mathfrak{H}$ to $\mathfrak{K}$ is a $CP$-semigroup of
$\mathfrak{B}
(\mathfrak{H} \oplus \mathfrak{K} )$.

Suppose $\gamma$ is a corner from $\alpha$ to $\beta$ and $\Theta$ is the
$CP$-semigroup of $\mathfrak{B}(\mathfrak{H} \oplus \mathfrak{K})$ defined
above. Suppose
$\Theta^\prime$ is a subordinate of $\Theta$ of the form
$$
\Theta_t^\prime \left[\begin{matrix} A&B
\\
C&D
\end{matrix} \right] =
\left[\begin{matrix} \alpha_t^\prime(A)&\gamma_t(B)
\\
\gamma_t^*(C)&\beta_t^\prime(D)
\end{matrix} \right]
$$
for $t \geq 0$ for $A,B,C$ and $D$ as stated above.  We say $\gamma$ is
\empha{maximal} if for every subordinate $\Theta^\prime$ of the above form we have
$\alpha^\prime = \alpha$.  We say $\gamma$ is \empha{hyper maximal} if for every
subordinate $\Theta^\prime$ of the above form we have $\alpha^\prime =
\alpha$ and $\beta^\prime = \beta$.
\end{definition}

We state Theorem~3.13 of \cite{BP} which shows how to determine when two
$CP$-semigroups dilate to cocycle  conjugate $E_0$-semigroups.

\begin{theorem}\label{thm:cocycle-conj-dil=hypermax-corner}%Theorem 1.6
Suppose $\alpha$ and $\beta$ are unital
$CP$-semigroups of $\mathfrak{B} (\mathfrak{H} )$ and $\mathfrak{B}
(\mathfrak{K} )$ and
$\alpha^d$ and $\beta^d$ are the minimal dilations of $\alpha$ and $\beta$
to $E_0$-semigroups.  Then $\alpha^d$ and $\beta^d$ are cocycle conjugate
if and only if there is a hyper maximal corner $\gamma$ from $\alpha$ to
$\beta$.
\end{theorem}

If $\alpha$ is a unital $CP$-semigroup and $\alpha^d$ is its minimal dilation
to an $E_0$-semigroup then the corners from $\alpha$ to $\alpha$ come from
contractive local cocycles.  The following theorem follows from Theorem
3.16 and Corollary 3.17 of \cite{BP}.

\begin{theorem}\label{thm:contractive-cocycles=eigencorners}%Theorem 1.7
Suppose $\alpha$ is a unital $CP$-semigroup of
$\mathfrak{B} (\mathfrak{H} )$ and $\alpha^d$ is its minimal dilation to an
$E_0$-semigroup $\alpha^d$ of $\mathfrak{B} (\mathfrak{H}_1)$. The relation
between
$\alpha$ and $\alpha^d$ is given by
$$
\alpha_t(A) = W^*\alpha_t^d(WAW^*)W
$$
for $A \in \mathfrak{B} (\mathfrak{H} )$ and $t \geq 0$ where $W$ is an isometry
from
$\mathfrak{H}$ to $\mathfrak{H}_1$ and $\alpha^d$ is minimal over the range of
$W$.
Suppose $\gamma$ is a corner from $\alpha$ to $\alpha$.  Then there is a
unique contractive local cocycle $C$ for $\alpha^d$ so that
$$
\gamma_t(A) = W^*C(t)\alpha_t^d(WAW^*)W
$$
for all $A \in \mathfrak{B} (\mathfrak{H} )$ and $t \geq 0$.  Conversely, if $C$
is a
contractive local cocycle for $\alpha^d$ then $\gamma$ given above is
corner from $\alpha$ to $\alpha $.

Furthermore, $C(t)$ is an isometry for all $t \geq 0$ if and only if
$\gamma$ is maximal and $C(t)$ is unitary for all $t \geq 0$ if and only if
$\gamma$ is hyper maximal.
\end{theorem}

Also in \cite{BP} (Theorem 3.16) there is a similar theorem for matrices of
corners.

\begin{theorem}\label{thm:matrix-contractive-cocycles=eigencorners}%Thm 1.8
Suppose $\alpha$ is a unital
$CP$-semigroup of
$\mathfrak{B} (\mathfrak{H} )$ and $\alpha^d$ is its minimal dilation to an
$E_0$-semigroup $\alpha^d$ of $\mathfrak{B} (\mathfrak{H}_1)$. The relation
between
$\alpha$ and $\alpha^d$ is given by
$$
\alpha_t(A) = W^*\alpha_t^d(WAW^*)W
$$
for $A \in \mathfrak{B} (\mathfrak{H} )$ and $t \geq 0$ where $W$ is an isometry
from
$\mathfrak{H}$ to $\mathfrak{H}_1$ and $\alpha^d$ is minimal over the range of
$W$.

Suppose $n$ is a positive integer and $\Theta$ is positive $(n \times
n)$-matrix of corners from $\alpha$ to $\alpha $.  Then there is a unique
positive $(n \times n)$-matrix $C$ of contractive local cocycles $C_{ij}$
for $\alpha^d$ for $i,j = 1,\dots ,n$ so that
$$
\theta_t^{(ij)}(A) = W^*C_{ij}(t)\alpha_t^d(WAW^*)W
$$
for all $A \in \mathfrak{B} (\mathfrak{H} )$ and $t \geq 0$.  Conversely, if $C$
is a
positive $(n \times n)$-matrix of contractive local cocycles for $\alpha^d$
then the matrix $\Theta$ whose coefficients $\theta^{(ij)}$ are given above
is a positive $(n \times n)$-matrix of corners from $\alpha$ to $\alpha $.
\end{theorem}

Next we define $CP$-flows.  We believe these are the simplest objects which
can be dilated to produce all spatial $E_0$-semigroups. $CP$-flows are
studied extensively in \cite{BP}.

\begin{definition}\label{def:CP-flow}%Definition 1.9
Suppose $\mathfrak{K}$ is a separable
Hilbert space
and $\mathfrak{H} = \mathfrak{K} \otimes L^2(0,\infty )$ and $U(t)$ is right
translation of $\mathfrak{H}$ by $t \geq 0$. Specifically, we may realize
$\mathfrak{H}$ as the space of $\mathfrak{K}$-valued Lebesgue measurable
functions with inner product
$$
(f,g) = \int_0^\infty (f(x),g(x)) \,dx
$$
for $f,g \in \mathfrak{H} $.  The action of $U(t)$ on an element $f \in
\mathfrak{H}$
is given by $(U(t)f)(x) = f(x - t)$ for $x \in [t,\infty )$ and $(U(t)f)(x)
= 0$ for $x \in [0,t)$.  A semigroup $\alpha$ is a
\empha{$CP$-flow over $\mathfrak{K}$}
if $\alpha$ is a $CP$-semigroup of $\mathfrak{B} (\mathfrak{H} )$ which is
intertwined by the translation semigroup $U(t)$, i.e., $U(t)A =
\alpha_t(A)U(t)$ for all $A \in \mathfrak{B} (\mathfrak{H} )$ and $t \geq 0$.
\end{definition}

In Theorem 4.0A of \cite{BP} it is shown that every spatial $E_0$-semigroup
is cocycle conjugate to an $E_0$-semigroup which is also a $CP$-flow.

We introduce notation for describing $CP$-flows.  Let $\mathfrak{H} =
\mathfrak{K}
\otimes L^2(0,\infty )$ and $U(t)$ be translation by $t.  $ Let
$$
E(t) = I - U(t)U(t)^*\qquad \text{and} \qquad E(a,b) = U(b)U(b)^*- U(a)U(a)
^*
$$
for $t \in [0,\infty )$ and $0 \leq a < b < \infty $.  Let $d = d/dx$ be
the differential operator of differentiation with the boundary condition
$f(0) = 0$. Precisely, the domain $\mathfrak{D} (d)$ is all $f \in \mathfrak{H}$
of
the form
$$
f(x) = \int_0^x g(t) \,dt
$$
with $g \in \mathfrak{H} $.  The hermitian adjoint $d^*$ is $-d/dx$  with no
boundary condition at $x = 0$.  Precisely, the domain $\mathfrak{D} (d^*)$
consists of the linear span of $\mathfrak{D} (d)$ and the functions $g(x) =
e^{-x}k$ with $k \in \mathfrak{K} $.

Suppose $\alpha$ is a $CP$-flow over $\mathfrak{K}$ and $A \in \mathfrak{B}
(\mathfrak{H}
)$. Then, for $t > 0$, one computes
$$
\alpha_t(A) = U(t)AU(t)^* + E(t)\alpha_t(A)E(t) = U(t)AU(t)^* + B
$$
for all $t \geq 0$.  Then $B$ commutes with $E(s)$ for all $s$ so $B$ is of
the form
$$
(Bf)(x) = b(x)f(x)
$$
and for $t > x \geq 0,\medspace b(x) \in \mathfrak{B} (\mathfrak{K} )$ depends
$\sigma$-strongly on $A$.  We define the boundary representation, $\pi_0$.
Let $\delta$ be the generator of $\alpha$. Then for $A \in \mathfrak{D} (\delta
)$ we have $A\mathfrak{D} (d) \subset \mathfrak{D} (d)$ and $A\mathfrak{D} (d^*)
\subset
\mathfrak{D} (d^*)$ so $A$ acts on $\mathfrak{D} (d^*)$ mod $\mathfrak{D} (d) =
\mathfrak{K}
$.  We call this mapping from $\mathfrak{D} (\delta )$ into $\mathfrak{B}
(\mathfrak{K} )$
the \empha{boundary representation} $\pi_0$.  Note $\pi_0$ tells you what flows in
from the origin.  The boundary representation need not be $\sigma$-weakly
continuous and even when it is it may not tell the whole story.  If $\pi$
is a $\sigma$-weakly continuous completely positive contraction of $\mathfrak{B}
(\mathfrak{K} \otimes L^2(0,\infty ))$ into $\mathfrak{B} (\mathfrak{K} )$ then
there is a
minimal $CP$-flow with that boundary representation and if that flow is
unital then the $E_0$-semigroup induced by the flow is completely spatial
(type~I$_n)$ where $n$ is the rank of $\pi $.

We now define the generalized boundary representation.  The resolvent
$R_\alpha$ for $\alpha$ is given by
$$
R_\alpha (A) = \int_0^\infty e^{-t}\alpha_t(A) \,dt
$$

Next we introduce some notation.  If $\phi$ is a $\sigma$-weakly continuous
mapping from $\mathfrak{B} (\mathfrak{H} )$ to $\mathfrak{B} (\mathfrak{K} )$ we
define $\hat
\phi$ is the predual map from $\mathfrak{B} (\mathfrak{K} )_*$ to $\mathfrak{B}
(\mathfrak{H}
)_*$ so we have $\rho (\phi (A)) = (\hat \phi\rho )(A)$ for all $A \in
\mathfrak{B} (\mathfrak{H} )$ and $\rho \in \mathfrak{B} (\mathfrak{K} )_*$. We
define the
mapping $\Gamma$ as
$$
\Gamma (A) = \int_0^\infty e^{-t} U(t)AU(t)^*dt
$$
for $A \in \mathfrak{B} (\mathfrak{H} )$.  Note $R_\alpha - \Gamma$ is
completely
positive which we denote by writing $R_\alpha - \Gamma \geq 0$.  Note
$\Gamma$ is the resolvent of a $CP$-flow with boundary representation
$\pi_0 = 0$.

We need one more bit of notation.  We define $\Lambda : \mathfrak{B}
(\mathfrak{K} )
\rightarrow \mathfrak{B} (\mathfrak{H} )$ for $A \in \mathfrak{B} (\mathfrak{K}
)$ we define
$\Lambda (A)$ by
$$
(\Lambda (A)f) = e^{-x}Af(x)
$$
We define $\Lambda = \Lambda (I)$.  Note $\Gamma (I) = I - \Lambda $.

Now we present the main formula.
$$
\hat R_\alpha (\rho ) = \hat \Gamma (\omega (\hat \Lambda\rho ) + \rho )
$$
for $\rho \in \mathfrak{B} (\mathfrak{H} )_*$ and $\eta \rightarrow \omega (\eta
)$
is the boundary weight map and $\omega (\eta )$ is the boundary weight
associated with $\eta $.  A boundary weight is a particular example of a
$T$-weight which we define presently.

\begin{definition}%{Definition 1.10}
Suppose $T \in \mathfrak{B} (\mathfrak{H}
)$ is a
positive strictly contractive operator (i.e. $ 0 \leq T \leq I$ and $\Vert
Tf \Vert < 1$ for $\Vert f \Vert \leq 1$ so one is not an eigenvalue for
$T)$.  We denote by $\mathfrak{A} (\mathfrak{H} ,T)$ the linear space
$$
\mathfrak{A} (\mathfrak{H} ,T) = (I - T)^{\tfrac{1}{2}}\mathfrak{B}
(\mathfrak{H} )(I - T)^{\tfrac{1}{2}}
$$
and by $\mathfrak{A} (\mathfrak{H} ,T)_*$ the linear functionals $\rho$ on
$\mathfrak{A} (
\mathfrak{H} ,T)$ of the form
$$
\rho ((I - T)^{\tfrac{1}{2}} A(I - T)^{\tfrac{1}{2}} ) = \eta (A)
$$
for $A \in \mathfrak{B} (\mathfrak{H} )$ with $\eta \in \mathfrak{B}
(\mathfrak{H} )_*$. We
call such functionals $T$-weights.  The $T$-norm of a $T$-weight $\rho$
denoted $\Vert\rho\Vert_T$ is the norm of $\eta$.  If $\rho$ is a
$T$-weight and $\Vert \rho\Vert_T \leq 1$ we say $\rho$ is $T$-contractive.
\end{definition}

Suppose $T \in \mathfrak{B}(\mathfrak{H})$ is a positive strictly contractive
operator and $P(\lambda )$ is the spectral resolution of $T$ so
$$
T = \int_0^1 \lambda dP(\lambda ).
$$
If $\rho \in \mathfrak{A} (\mathfrak{H} ,T)_*$ then $\rho$ restricted to $
P(\lambda
)\mathfrak{B} (\mathfrak{H} )P(\lambda )$ is normal for all $\lambda > 0$.

Consider now the case when $T_1 \geq T_2 \geq 0$ and $T_1$ is strictly
contractive so $T_2$ is strictly contractive.  Define $S$ on the range
$\sqrt{I - T_2}$ by the relation
$$
S(I - T_2)^{\tfrac{1}{2}} f = (I - T_1)^{\tfrac{1}{2}} f
$$
for $f \in \mathfrak{H}$.  Note
$$
\Vert S(I - T_2)^{\tfrac{1}{2}} f\Vert^2 = (f,(I - T_1)f) \leq (f,(I - T_2)f) =
\Vert (I - T_2)^{\tfrac{1}{2}} f\Vert^2
$$
for $f \in \mathfrak{H}$.  Then $S$ is a contractive map on the range of
$\sqrt{I - T_2}$ which is dense in $\mathfrak{H}$  so $S$ has a unique bounded
extension to a contraction defined on all of $\mathfrak{H}$.  We also denote
this operator by $S$.  We note $S$ is a contraction which satisfies the
operator equation
$$
S(I - T_2)^{\tfrac{1}{2}} = (I - T_1)^{\tfrac{1}{2}} \quad\; \text{so} \quad (I
- T_2)^{\tfrac{1}{2}} S^*AS(I - T_2)^{\tfrac{1}{2}} = (I - T_1)^{\tfrac{1}{2}}
A( - T_1)^{\tfrac{1}{2}}
$$
for $A \in \mathfrak{B} (\mathfrak{H} )$ so it follows that $\mathfrak{A}
(\mathfrak{H} ,T_1)
\subset \mathfrak{A} (\mathfrak{H} ,T_2)$.  We show that $\mathfrak{A}
(\mathfrak{H} ,T_2)_*
\subset \mathfrak{A} (\mathfrak{H} ,T_1)_*$.  Suppose $\rho \in \mathfrak{A}
(\mathfrak{H}
,T_2)_*$ which means
$$
\rho ((I - T_2)^{\tfrac{1}{2}} A(I - T_2)^{\tfrac{1}{2}} ) =\eta (A)
$$
for all $A \in \mathfrak{B} (\mathfrak{H} )$ where $\eta \in \mathfrak{B}
(\mathfrak{H} )_*$.
Then we have
$$
\rho ((I - T_1)^{\tfrac{1}{2}} A(I - T_1)^{\tfrac{1}{2}} ) = \rho ((I -
T_2)^{\tfrac{1}{2}}
S^*AS(I - T_2)^{\tfrac{1}{2}} ) = \eta (S^*AS)
$$
for $A \in \mathfrak{B} (\mathfrak{H} )$.  So we see $\rho \in \mathfrak{A}
(\mathfrak{H}
,T_1)$ and since $S$ is a contraction we have $\Vert\rho\Vert_{T_1} \leq
\Vert\rho\Vert_{T_2}$.

Note a $0$-weight is just a normal functional.

We caution the reader that the $T$-weights we consider are not normal
weights.  A normal weight on a von Neumann algebra has the property that if
$0 \leq A_1 \leq A_2 \leq \cdots $ is a increasing sequence of operators
which converge strongly to $A$  then
$\omega (A)$ is the limit of the $\omega (A_k)$, where we allow $+ \infty$
as a possible limit.  Let $\mathfrak{H} = L^2(0,\infty )$ and let $\omega$ be
the $\Lambda$-weight given by
$$
\omega ((I - \Lambda )^{\tfrac{1}{2}} A(I - \Lambda )^{\tfrac{1}{2}} ) = (h,Ah)
$$
for $A \in \mathfrak{B}(\mathfrak{H})$ where $h(x) = x^{-{\tfrac{1}{2}} s}(1 -
e^{-x})^{\tfrac{1}{2}}$ for $s \in (1,2)$.  Let $\mathfrak{M}_n$ be the
functions $g$ in $\mathfrak{H}$ with support in $[1/n,\infty )$ and
$$
\int_{1/n}^\infty x^{-{\tfrac{1}{2}} s}(1 - e^{-x})^{\tfrac{1}{2}} g(x) \,dx = 0
$$
and let $P_n$ be the orthogonal projection onto $\mathfrak{M}_n$ for $n =
1,2,\cdots $.  Note $\omega (P_n) = 0$ for $n = 1,2,\cdots $ and $P_n
\rightarrow I$ as $n \rightarrow \infty$ but it is not true that $\omega
(I) = 0$.  If we were to assign $\omega (I)$ a value it is $+\infty$ since
$\omega$ is positive and unbounded.  Although $T$-weights are not in
general normal weights we do not think of them as pathological like non-normal
bounded functionals since $T$-weights are normal when scaled down by
$\sqrt{I - T}$.

In the particular case when $\mathfrak{H} = \mathfrak{K} \otimes L^2(0,\infty )$
when
we speak of \empha{boundary weights} we mean the following.  Let $\Lambda$ be the operator corresponding to multiplication by $e^{-x}$.  Then the boundary
algebra $\mathfrak{A} (\mathfrak{H} )$ is
$$
\mathfrak{A} (\mathfrak{H} ) = \mathfrak{A} (\mathfrak{H} ,\Lambda) = (I -
\Lambda )^{\tfrac{1}{2}}
\mathfrak{B}(\mathfrak{H} )(I - \Lambda )^{\tfrac{1}{2}}
$$
and the boundary weights denoted by $\mathfrak{A} (\mathfrak{H} )_*$ are
$$
\mathfrak{A} (\mathfrak{H} )_* = \mathfrak{A} (\mathfrak{H} ,\Lambda)_*.
$$
If $\omega$ is a boundary weight we say $\omega$ is weight contractive if
$\Vert \omega\Vert_\Lambda \leq 1$ and if $\omega$ is a positive boundary
weight we say $\omega$ is normalized if $\Vert\omega\Vert_\Lambda = 1$.  If
$\omega$ is a boundary weight and we say $\omega$ is \empha{bounded} we mean
$\omega$ is bounded as a functional on $\mathfrak{B} (\mathfrak{H} )$ (i.e.
$\vert\omega (A)\vert \leq \Vert A\Vert$ for all $A \in \mathfrak{A}
(\mathfrak{H}
))$.

The mapping $\rho \rightarrow \omega(\rho)$ defined for $\rho \in \mathfrak{B}
(\mathfrak{K})_*$ is a \empha{boundary weight map} if this mapping is a linear mapping of $\mathfrak{B}(\mathfrak{K})_*$ into boundary weights on $\mathfrak{A}(\mathfrak{H})$ and
this mapping is completely bounded with the norm on
$\mathfrak{B}(\mathfrak{K})_*$
being the usual norm and the norm on the boundary weights being the
boundary weight norm.  A boundary weight map is positive if it is
completely positive.  A boundary weight map $\omega$ is unital if
$\omega(\rho)(I - \Lambda) = \rho(I)$ for all $\rho \in \mathfrak{B}(
\mathfrak{K})_*$.

Maintaining the notation of the above definition we observe that $U(t)AU(t)^*
\in \mathfrak{A}(\mathfrak{H})$ for all $A \in \mathfrak{B}(\mathfrak{H})$ and
$t > 0$. Recall
the mapping $\Gamma$ defined above.  Since $\Gamma$ is completely positive
and $\Gamma(I) = I - \Lambda$, so $\Gamma(I) \in \mathfrak{A}(\mathfrak{H})$, it
follows that $\Gamma(A) \in \mathfrak{A}(\mathfrak{H})$ for all $A \in
\mathfrak{B}(\mathfrak{H})$.  For more details see the discussion after
Definition 4.16 of \cite{BP}.

Every $CP$-flow is given by a boundary weight map $\rho \rightarrow \omega
(\rho )$.  As we have mentioned the map is completely positive.  There is a
further complicated positivity condition.  The condition says if you
construct an approximation to the boundary representation $\pi_t$, then
$\pi _t$ is completely positive.

We describe the connection between boundary weight and boundary
representation.  One can construct a boundary weight map so that the
boundary representation is a given $\sigma$-weakly continuous completely
positive contraction of $\mathfrak{B} (\mathfrak{H} )$ into $\mathfrak{B}
(\mathfrak{K} )$.
Suppose $\pi$ is a $\sigma$-weakly continuous completely positive
contraction of $\mathfrak{B} (\mathfrak{H} )$ into $\mathfrak{B} (\mathfrak{K}
)$. Let
$$
\omega = \hat \pi + \hat \pi\hat \Lambda\hat \pi + \hat \pi\hat \Lambda\hat
\pi\hat \Lambda\hat \pi + \hat \pi\hat \Lambda\hat \pi\hat \Lambda\hat
\pi\hat \Lambda\hat \pi + \dots
$$
This converges as a weight (i.e. the above series converges on the boundary
algebra $\mathfrak{A}(\mathfrak{H})$)  and this is the boundary weight of a
$CP$-flow.  We call this the minimal $CP$-flow derived from $\pi $.
Formally $\omega = \hat \pi (I - \hat \Lambda\hat \pi )^{-1}$ and solving
for $\pi$ we have
$$
\hat \pi = \omega (I + \hat \Lambda\omega )^{-1}
$$

If a boundary weight associated with a $CP$-flow is bounded the boundary
representation is well defined as stated in the next theorem (see theorem
4.27 of \cite{BP}).

\begin{theorem}\label{thm:bounded-bdry-wgt-->bdry-rep}%{Theorem 1.11}
Suppose $\alpha$ is a $CP$-flow over $\mathfrak{K}$ and
$\rho \rightarrow \omega (\rho )$ is the associated boundary weight map.
Suppose $\Vert\omega (\rho )\Vert < \infty$ for $\rho \in \mathfrak{B}
(\mathfrak{K}
)_*$ so $\omega (\rho ) \in \mathfrak{B} (\mathfrak{H} )_*$ for all $\rho
\in
\mathfrak{B} (\mathfrak{K} )_*$.  Then the mapping $\rho \rightarrow \rho + \hat
\Lambda\omega (\rho )$ is invertible $i.e. (I + \hat \Lambda\omega )^{-1}$
exists and $\hat \pi$ given by
$$
\hat \pi = \omega (I + \hat \Lambda\omega )^{-1}
$$
is a completely positive contraction from $\mathfrak{B} (\mathfrak{K} )_*$ to
$\mathfrak{B} (\mathfrak{H} )_*$.  There is a unique $CP$-flow derived from
$\pi$ and its boundary weight map is given by
$$
\omega = \hat \pi + \hat \pi\hat \Lambda\hat \pi + \hat \pi\hat \Lambda\hat
\pi\hat \Lambda\hat \pi + \hat \pi\hat \Lambda\hat \pi\hat \Lambda\hat
\pi\hat \Lambda\hat \pi + \dots
$$
So when $\omega (\rho )$ is bounded for all $\rho \in \mathfrak{B} (\mathfrak{K}
)_*
$ we have
$$
\omega = \hat \pi (I - \hat \Lambda\hat \pi )^{-1}\qquad \text{and} \qquad
\hat \pi = \omega (I + \hat \Lambda\omega )^{-1}
$$
\end{theorem}

Now we introduce a bit of notation.  Suppose $\omega$ is a boundary weight
and $t > 0$.  We denote by $\omega\vert_t$ the functional given
by $\omega\vert_t(A)
= \omega (E(t,\infty )AE(t,\infty ))$ for $A \in \mathfrak{B} (\mathfrak{H} )$.
Note
$\omega\vert_t(\rho ) \in \mathfrak{B} (\mathfrak{H} )_*$, i.e.
$\omega\vert_t(\rho
)$ is a bounded $\sigma$-weakly continuous functional.  We use the same
notation for operators.  If $A \in \mathfrak{B} (\mathfrak{H} )$ and $t > 0$
then we
denote $A\vert_t$ the operator $A\vert_t = E(t,\infty )AE(t,\infty )$.
Note for $\omega$ a boundary weight and $A \in \mathfrak{B} (\mathfrak{H} )$
then
$\omega\vert_t(A) = \omega (A\vert_t)$.

From Theorems 4.23 and 4.27 and Lemma 4.34 of \cite{BP} we have the following
theorem.

\begin{theorem}\label{thm:restricted-weights}%{Theorem 1.12}
Suppose $\rho \rightarrow \omega (\rho )$ is the
boundary weight map of a $CP$-flow over $\mathfrak{K} $.  Then for each $t > 0$
we have $\rho \rightarrow \omega\vert_t(\rho )$ is the boundary weight map
of a $CP$-flow over $\mathfrak{K} $.  Suppose $\rho \rightarrow \omega (\rho )$
is a completely positive mapping of $\mathfrak{B} (\mathfrak{K} )$ into boundary
weights on $\mathfrak{B} (\mathfrak{H} )$ satisfying $\omega (\rho )(I - \Lambda
)
\leq \rho (I)$ for $\rho$ positive.  Suppose
$$
\hat \pi_t^\# = \omega\vert_t(I + \hat \Lambda\omega\vert_t)^{-1}
$$
is a completely positive contraction of $\mathfrak{B} (\mathfrak{K} )_*$ into
$\mathfrak{B} (\mathfrak{H} )_*$ for each $t > 0$.  Then $\rho \rightarrow
\omega (\rho )$
is the boundary weight map of a $CP$-flow over $\mathfrak{K}$.

Furthermore, the mapping $\hat \pi_t^\#$ defined above has the property
that if $\phi_t(A) = \pi_t^\#(E(s,\infty)AE(s,\infty)$ for $0 < t \leq s <
\infty$ and $A \in \mathfrak{B}(\mathfrak{H})$ then $\phi_t$ is increasing in
$t$ in
the sense complete positivity (i.e., the mapping $A \rightarrow \phi_t(A) -
\phi_r(A)$ for $A \in \mathfrak{B}(\mathfrak{H})$ and $0 < t < r \leq s$ is
completely positive.
\end{theorem}

\begin{definition}% {Definition 1.13}
If $\rho \rightarrow \omega (\rho )$ is a
mapping of $\mathfrak{B} (\mathfrak{K} )_*$ into boundary weights on
$\mathfrak{B} (\mathfrak{H} )$ so that $\hat \pi_t^\#$ defined above is
completely positive for each $t > 0$ we say this map is \empha{$q$-positive}.  The
family $\pi_t^\#$ of completely positive $\sigma$-weakly continuous contractions
of $\mathfrak{B} (\mathfrak{H} )$ into
$\mathfrak{B} (\mathfrak{K} )$ is called the \empha{generalized boundary
representation}.
\end{definition}

We remark that in checking that the $\pi_t^\#$ are completely positive it
is only necessary to check for small $t$.  If the mapping $\pi_t^\#$ is
completely positive then $\pi_s^\#$ is completely positive for all $s \geq
t$.   Next we give the order relation for the generalized boundary
representation(see Theorem 4.20 of \cite{BP}).

\begin{theorem}%{Theorem 1.14}
If $\alpha$ and $\beta$ are $CP$-flows over
$\mathfrak{K}$ then $\beta$ is a subordinate of $\alpha$ ($\alpha \geq \beta$)
if and only if $\pi_t^\# \geq \phi_t^\#$ for all $t > 0$ where $\pi_t^\#$
and $\phi_t^\#$ are the generalized boundary representations of $\alpha$
and $\beta $.   Also we have if $\pi_t^\# \geq \phi_t^\#$ then $\pi_s^\#
\geq \phi_s^\#$ for all $s \geq t$ so one only has to check for a sequence
$\{t_n\}$ tending to zero.
\end{theorem}

In Theorem~\ref{thm:bounded-bdry-wgt-->bdry-rep} we used the phrase
$\alpha$ is derived from $\pi$.  The next
theorem (see Theorem 4.24 of \cite{BP})and definition will make this more
precise.  We need a bit of notation which is given in the next definition.

\begin{definition}%{Definition 1.15}
Let $Q_\subo$ be the map from $\mathfrak{K}$
to $\mathfrak{H}$ given by $(Q_\subo k)(x) = e^{-\sqrt{x}} k$.  And let $\Phi$
be the
mapping of $\mathfrak{B} (\mathfrak{K} )_*$ into $\mathfrak{B} (\mathfrak{H}
)_*$ given by
$\Phi (\rho )(A) = \rho (Q_\subo ^*AQ_\subo )$ for all $A \in \mathfrak{B}
(\mathfrak{H}
)$.
\end{definition}

Note that $\Phi (\rho )(U(t)AU(t)^*) = e^{-t}\Phi (\rho )(A)$, $\Phi (\rho
)(\Gamma (A)) = {\tfrac{1}{2}}\Phi (\rho )(A)$ for $t \geq  0$ and $\Phi (\rho
)(\Lambda (C)) = {\tfrac{1}{2}}\rho (C)$ for all $\rho \in \mathfrak{B}
(\mathfrak{K} )_*$,
$A \in \mathfrak{B} (\mathfrak{H} )$ and $C \in \mathfrak{B} (\mathfrak{K} )$.

\begin{theorem}\label{thm:boundary-weight-conditions}%{Theorem 1.16}
Suppose $\rho \rightarrow \omega (\rho )$ defines
a $CP$-flow over $\mathfrak{K}$ as described in Definition~\ref{def:CP-flow} and
$\delta$ is the
generator of $\alpha$ (i.e., $\delta$ is the derivative of $\alpha_t$ at
$t = 0$).  Suppose $\pi$ is a completely positive normal
contraction of $\mathfrak{B} (\mathfrak{H} )$ into $\mathfrak{B} (\mathfrak{K}
)$. Then the
following are equivalent.
\begin{enumerate}
\item $\Phi (\rho ) \in \mathfrak{D} (\hat \delta )$ and $\hat \delta (\Phi
(\rho )) = \hat \pi (\rho ) - \Phi (\rho )$ for each $\rho \in \mathfrak{B}
(\mathfrak{K} )_*$.

\item $\omega (\rho - \hat \Lambda (\hat \pi (\rho ))) = \hat \pi
(\rho )$ for all $\rho \in \mathfrak{B} (\mathfrak{K} )_*$.

\item $\pi (A) = \pi_\subo (A)$ for all $A \in \mathfrak{D} (\delta )$ where
$\pi _\subo $ is the boundary representation of $\alpha$.
\end{enumerate}
\end{theorem}

\begin{definition}\label{def:derived-CP-flow}%{Definition 1.17 }
We say a $CP$-flow $\alpha$ over $\mathfrak{K}$
is derived from the completely positive normal contraction $\pi$ of
$\mathfrak{B} (\mathfrak{H} )$ into $\mathfrak{B} (\mathfrak{K} )$ if it
satisfies one and, therefore,
all the conditions of Theorem~\ref{thm:boundary-weight-conditions}.
\end{definition}

As mentioned earlier for each such $\pi$ there is a $CP$-flow $\alpha$
derived from $\pi$ and the next theorem (see Theorem 4.26 of \cite{BP})
gives a condition for uniqueness.

\begin{theorem}\label{thm:bdry-wght-when-rep-normal}%{Theorem 1.18}
Suppose $\pi$ is a completely positive
$\sigma$-weakly continuous linear contraction of $\mathfrak{B} (\mathfrak{H} )$
into
$\mathfrak{B} (\mathfrak{K} )$.  Then for each $\rho \in \mathfrak{B}
(\mathfrak{K} )_*$ the
sum
$$
\omega (\rho ) = \hat \pi (\rho ) + \hat \pi (\hat \Lambda (\hat \pi (\rho
))) + \hat \pi (\hat \Lambda (\hat \pi (\hat \Lambda (\hat \pi (\rho )))))
+ \cdots
$$
converges as a weight on $\mathfrak{A} (\mathfrak{H})$ and the mapping $\rho
\rightarrow \omega (\rho )$ is the boundary weight map of a $CP$-flow
$\alpha$ which is derived from $\pi$.  Furthermore, this $\alpha$ is the
minimal $CP$-flow derived from $\pi$ in that if $\rho \rightarrow \eta
(\rho )$ is the boundary weight map of a second $CP$-semigroup derived from
$\pi$ then $\omega (\rho ) \leq \eta (\rho )$ for all positive $\rho \in
\mathfrak{B} (\mathfrak{K} )_*$.  Moreover, if $(\pi\circ\Lambda )^n(I)$
$\rightarrow 0$ weakly as $n \rightarrow \infty$ then $\alpha$ defined
above is unique (i.e.  $\alpha$ is the only $CP$-flow derived from $\pi
)$.
\end{theorem}

We remark that we believe that this theorem can be strengthened with the
stronger conclusion being that $\alpha$ is a flow subordinate of any
$CP$-flow derived from $\pi$.  In the examples we construct in this paper
we will show that the stronger result holds.

So far most of the results in this section are proved in \cite{BP}.  The
next two theorems are new.

\begin{theorem}\label{thm:derivation-and-subordination}  %{Theorem 1.19}
Suppose $\alpha$ is a $CP$-flow over $\mathfrak{K}$
derived from $\pi$ as described in Definition~\ref{def:derived-CP-flow} and
$\beta$ is $CP$-flow
subordinate to $\alpha $, so the mapping $A \rightarrow \alpha_t(A) -
\beta_t(A)$ for $A \in \mathfrak{B} (\mathfrak{H} )$ is completely positive for
all
$t \geq 0$.  Then there is a unique completely positive normal contraction
$\phi$ of $\mathfrak{B} (\mathfrak{H} )$ into $\mathfrak{B} (\mathfrak{K} )$
which is
subordinate to $\pi$ so that $\beta$ is derived from $\phi$.
\end{theorem}

\begin{proof}  Assume the hypothesis of the theorem and suppose
$\delta_\alpha$ and $\delta_\beta$ are the generators of $\alpha$ and
$\beta$, respectively.  Let $\gamma_t(A) = U(t)AU(t)^*$ for $t \geq 0$ and
$A \in \mathfrak{B} (\mathfrak{H} )$.  Since $\beta$ is a subordinate of
$\alpha$ and
$\beta$ is intertwined by $U(t)$ we have the maps $t \rightarrow
\alpha_t(A) - \beta_t(A)$ and $t \rightarrow$ $\beta_t(A) - \gamma_t(A)$
and $A \in \mathfrak{B} (\mathfrak{H} )$ are completely positive for all
$t \geq 0$.
Suppose $\rho \in \mathfrak{B} (\mathfrak{K} )_*$ and $\rho \geq 0$.  Then we
have
\begin{align*}
\vartheta_t &= t^{-1}(\hat \alpha_t(\Phi (\rho )) - \Phi (\rho )) + \Phi
(\rho ) \geq t^{-1}(\hat \beta_t(\Phi (\rho )) - \Phi (\rho )) + \Phi
(\rho )
\\
&= \nu_t  \geq t^{-1}(\hat \gamma_t(\Phi (\rho )) - \Phi (\rho )) + \Phi
(\rho ) = t^{-1}(e^{-t} - 1 + t)\Phi (\rho )
\end{align*}
for $t > 0$ where the two equal signs are definitions of $\vartheta_t$ and
$\nu _t$.  Since $\alpha$ is derived from $\pi$ we have
$$
\vartheta_t = t^{-1}(\hat \alpha_t(\Phi (\rho )) - \Phi (\rho )) + \Phi
(\rho ) \rightarrow \hat \delta (\Phi (\rho )) + \Phi (\rho ) = \hat \pi
(\rho ) = \vartheta_\subo
$$
as $t \rightarrow 0^+$ and the convergence is in norm.  Since $\vartheta_\subo
= \hat \pi (\rho ) \in \mathfrak{B} (\mathfrak{H} )_*$ there is a positive trace
class operator $\Omega_\subo $ so that
$$
\vartheta_\subo (A) = \hat \pi (\rho )(A) = tr(A\Omega_\subo )
$$
for $A \in \mathfrak{B} (\mathfrak{H} )$ and for every $\rho_1 \in \mathfrak{B}
(\mathfrak{H} )
_*$ with $\vartheta_\subo  \geq \rho_1 \geq 0$ there is an $X \in \mathfrak{B}
(\mathfrak{H} )$ with $0 \leq X \leq I$ so that
$$
\rho_1(A) = tr(A\Omega_\subo ^{\tfrac{1}{2}} X\Omega_\subo ^{\tfrac{1}{2}} )
$$
for $A \in \mathfrak{B} (\mathfrak{H} )$ and conversely if $X \in \mathfrak{B}
(\mathfrak{H} )$
and $0 \leq X \leq I$ then $\rho_1$ defined above is in $\mathfrak{B}
(\mathfrak{H}
)_*$ and $0 \leq \rho_1 \leq \vartheta_\subo $.  (If we require the null space
of $X$ contains Range($\Omega_\subo )^\perp$ then $X$ is uniquely determined by
$\rho _1.)$ Suppose $t > 0$ and $\Omega_t$ is the unique positive trace
class operator so that
$$
\nu_t(A) = (t^{-1}(\hat \alpha_t(\Phi (\rho )) - \Phi (\rho )) + \Phi (\rho
))(A) = tr(A\Omega_t)
$$
for $A \in \mathfrak{B} (\mathfrak{H} )$.  From the inequality above we have
$\vartheta _t \geq \nu_t \geq 0$ so there is an operator $X_t \in \mathfrak{B}
(\mathfrak{H} )$ with $0 \leq X_t \leq I$ so that
$$
\nu_t(A) = (t^{-1}(\hat \beta_t(\Phi (\rho )) - \Phi (\rho )) + \Phi (\rho
))(A) = tr(A\Omega_t^{\tfrac{1}{2}} X_t\Omega_t^{\tfrac{1}{2}} )
$$
for $A \in \mathfrak{B} (\mathfrak{H} )$.  We will require the null space of
$X_t$
contains Range($\Omega_t)^\perp$ so $X_t$ is uniquely determined.  Now let
$$
\eta_t(A) = tr(A\Omega_\subo ^{\tfrac{1}{2}} X_t\Omega_\subo ^{\tfrac{1}{2}} )
$$
for $A \in \mathfrak{B} (\mathfrak{H} )$.  Now we have
$$
\Vert\nu_t - \eta_t\Vert = \sup\{ Re(tr(A(\Omega_t^{\tfrac{1}{2}}
X_t\Omega_t^{\tfrac{1}{2}} - \Omega_\subo ^{\tfrac{1}{2}} X_t\Omega_\subo
^{\tfrac{1}{2}}
)): A \in
\mathfrak{B} (\mathfrak{H} ),\medspace \Vert A\Vert \leq 1\} .
$$
We have $\vert tr(AB)\vert \leq  \Vert A\Vert_{HS} \Vert B\Vert _{HS}$ for
$A,B \in \mathfrak{B} (\mathfrak{H} )$ with $\Vert A\Vert_{HS} =
tr(A^*A)^{\tfrac{1}{2}}$
the Hilbert Schmidt norm we have for $A \in \mathfrak{B} (\mathfrak{H} )$ with
$\Vert A\Vert \leq 1$ that
\begin{align*}
\vert tr(A(\Omega_t^{\tfrac{1}{2}} X_t\Omega_t^{\tfrac{1}{2}} -
& \Omega_\subo ^{\tfrac{1}{2}}
X_t\Omega _\subo ^{\tfrac{1}{2}} )\vert
\\
&= \vert tr(A(\Omega_t^{\tfrac{1}{2}} X_t\Omega_t^{\tfrac{1}{2}} - \Omega
_t^{\tfrac{1}{2}}
X_t\Omega_\subo ^{\tfrac{1}{2}}  + \Omega_t^{\tfrac{1}{2}} X_t\Omega_\subo
^{\tfrac{1}{2}}
-
\Omega_\subo ^{\tfrac{1}{2}} X_t\Omega_\subo ^{\tfrac{1}{2}} )\vert
\\
&\leq \vert tr(A\Omega_t^{\tfrac{1}{2}} X_t(\Omega_t^{\tfrac{1}{2}} -
\Omega_\subo ^{\tfrac{1}{2}}
))\vert + \vert tr(X_t\Omega_\subo ^{\tfrac{1}{2}} A(\Omega_t^{\tfrac{1}{2}} -
\Omega_\subo ^{\tfrac{1}{2}} ))\vert
\\
&\leq (\Vert A\Omega_t^{\tfrac{1}{2}} X_t\Vert_{HS} + \Vert
A\Omega_\subo ^{\tfrac{1}{2}}
X_t \Vert_{HS}) \Vert\Omega_t^{\tfrac{1}{2}} - \Omega_\subo
^{\tfrac{1}{2}}\Vert_{HS}
\\
&\leq (tr(\Omega_t)^{\tfrac{1}{2}} + tr(\Omega_\subo )^{\tfrac{1}{2}}
)\Vert\Omega_t^{\tfrac{1}{2}} - \Omega_\subo ^{\tfrac{1}{2}}\Vert_{HS}
\end{align*}
and, hence, it follows that
$$
\Vert\nu_t - \eta_t\Vert \leq (\Vert\nu_t\Vert^{\tfrac{1}{2}} +
\Vert\vartheta_\subo \Vert ^{\tfrac{1}{2}} )\Vert\Omega_t^{\tfrac{1}{2}} -
\Omega_\subo ^{\tfrac{1}{2}}\Vert_{HS} \leq (\Vert
\vartheta_t\Vert^{\tfrac{1}{2}} +
\Vert\vartheta_\subo \Vert^{\tfrac{1}{2}} )\Vert\Omega_t^{\tfrac{1}{2}} -
\Omega_\subo ^{\tfrac{1}{2}}\Vert_{HS}
$$
Now if $UT$ is the polar decomposition of $\Omega_t^{\tfrac{1}{2}} -
\Omega_\subo ^{\tfrac{1}{2}}$ so $U^*(\Omega_t^{\tfrac{1}{2}} -
\Omega_\subo ^{\tfrac{1}{2}} ) =$
$\vert\Omega_t^{\tfrac{1}{2}} - \Omega_\subo ^{\tfrac{1}{2}}\vert =
((\Omega_t^{\tfrac{1}{2}} -
\Omega _\subo ^{\tfrac{1}{2}} )^2)^{\tfrac{1}{2}}$ we have
\begin{align*}
\Vert\vartheta_t - \vartheta_\subo \Vert &\geq \vert\vartheta_t(U^*) - \vartheta
_\subo (U^*)\vert = \vert tr(U^*(\Omega_t - \Omega_\subo ))\vert
\\
&= \vert tr(U^*((\Omega_t^{\tfrac{1}{2}} - \Omega_\subo ^{\tfrac{1}{2}}
)(\Omega_t^{\tfrac{1}{2}}
+ \Omega_\subo ^{\tfrac{1}{2}} ) + (\Omega_t^{\tfrac{1}{2}} + \Omega_\subo
^{\tfrac{1}{2}}
)(\Omega_t^{\tfrac{1}{2}} - \Omega_\subo ^{\tfrac{1}{2}} )))\vert
\\
&= tr(\vert\Omega_t^{\tfrac{1}{2}} - \Omega_\subo ^{\tfrac{1}{2}}\vert
(\Omega_t^{\tfrac{1}{2}} +
\Omega_\subo ^{\tfrac{1}{2}} )) \geq tr(\vert\Omega_t^{\tfrac{1}{2}} -
\Omega_\subo ^{\tfrac{1}{2}}\vert ^2)
\\
&= tr((\Omega_t^{\tfrac{1}{2}} - \Omega_\subo ^{\tfrac{1}{2}} )^2) =
\Vert\Omega_t^{\tfrac{1}{2}}
- \Omega_\subo ^{\tfrac{1}{2}}\Vert_{HS}^2.
\end{align*}
Hence, we have
$$
\Vert\nu_t - \eta_t\Vert \leq (\Vert\vartheta_t\Vert^{\tfrac{1}{2}} +
\Vert\vartheta _\subo \Vert^{\tfrac{1}{2}} )\Vert\vartheta_t -
\vartheta_\subo \Vert^{\tfrac{1}{2}}
$$
for all $t > 0$.  Note $\vartheta_\subo  \geq \eta_t \geq 0$ for each $t > 0$.
Note the set $S$ of $\eta \in \mathfrak{B} (\mathfrak{H} )_*$ with
$\vartheta_\subo
\geq
\eta \geq 0$ is compact.  This may be seen as follows.  For every
$\epsilon_1 > 0$ there is a finite rank $\xi \in \mathfrak{B}(\mathfrak{H})_*$
so
that $0 \leq \xi \leq \vartheta_\subo $ and $\Vert\xi - \vartheta_\subo \Vert <
\epsilon_1$.  Given $\epsilon > 0$ we can by choosing $\epsilon_1$ small
enough insure that for every $\eta \in S$ there is a positive $\eta^\prime
\leq \xi$ with $\Vert \eta^\prime - \xi \Vert < \epsilon$.  Hence, for
every $\epsilon > 0$ there is a finite dimensional compact subset of $S$
so that the $\epsilon$-neighborhoods of this set cover $S$, thus for every
$\epsilon > 0$ there is a cover of S with a finite numbers of open balls of
radius $\epsilon$ and we obtain that $S$ is compact.  Since $S$ is compact there is a
sequence $t_n \rightarrow 0^+$ as $n \rightarrow \infty$ so that
$\eta_{t_n}$ converges to a limit $\eta_\subo $ in norm as $n \rightarrow \infty
$.  Since $\Vert\vartheta_t - \vartheta_\subo \Vert \rightarrow 0$ as $t
\rightarrow 0^+$ it follows from the above estimate that $\nu_{t_n}
\rightarrow \eta_\subo $ as $n \rightarrow \infty$.  Hence, we have
$$
\Vert t_n^{-1}(\hat \beta_{t_n}(\Phi (\rho )) - \Phi (\rho )) - (\eta_\subo  -
\Phi (\rho ))\Vert \rightarrow 0
$$
as $n \rightarrow \infty$.  Now let
$$
\mu_n = \frac {1} {t_n} \int_0^{t_n} \hat \beta_s(\Phi (\rho )) \,ds
$$
for $n = 1,2,\cdots$.  Let $\delta_\beta$ be the generator of $\beta$.
Note $\mu _n \in \mathfrak{D} (\hat \delta_\beta )$ and
$$
\hat \delta_\beta (\mu_n) = t_n^{-1}(\hat \beta_{t_n}(\Phi (\rho )) - \Phi
(\rho )) \rightarrow \eta_\subo  - \Phi (\rho )
$$
as $n \rightarrow \infty$ where the convergence is in norm.  Since $\mu_n
\rightarrow \Phi (\rho )$ in norm as $n \rightarrow \infty$ it follows from
the fact that $\hat \delta_\beta$ is closed that $\Phi (\rho ) \in \mathfrak{D}
(\hat \delta_\beta )$ and $\hat \delta_\beta (\Phi (\rho )) = \eta_\subo  - \Phi
(\rho )$.  Since $\rho$ was an arbitrary positive element of $\mathfrak{B}
(\mathfrak{H} )_*$ it follows that $\Phi (\rho ) \in \mathfrak{D} (\hat \delta
_\beta
)$ and $\hat \delta_\beta (\Phi (\rho )) + \Phi (\rho ) = \hat \phi (\rho
)$ where this equation defines $\phi$.  Since, $\alpha \geq \beta \geq
\gamma$ in the sense of complete positivity it follows that $\pi \geq
\phi \geq 0$ is the sense of complete positivity.   From
Definition~\ref{def:derived-CP-flow} it
follows that $\beta$ is derived from $\phi$.  The uniqueness of $\phi$
follows from the defining equation for $\phi$.
\end{proof}

\begin{theorem}%{Theorem 1.20}
Suppose $\alpha$ is a $CP$-flow over
$\mathfrak{K}$ and
$\pi$ is a normal completely positive contraction of $\mathfrak{B} (\mathfrak{H}
)$
into $\mathfrak{B} (\mathfrak{K} )$ and suppose further that $\pi$ is unital so
$\pi
(I) = I$.  Suppose $\beta$ is a $CP$-flow over $\mathfrak{K}$ derived from $\pi$
and $\alpha \geq \beta$ (i.e. the mapping $A \rightarrow \alpha_t(A) -
\beta_t(A)$ for $A \in \mathfrak{B} (\mathfrak{H} )$ is completely positive for
all
$t \geq 0)$.  Then $\alpha$ is derived from $\pi $.
\end{theorem}
\begin{proof}  Assume the hypothesis and notation of the theorem.  Suppose
$\rho \in \mathfrak{B} (\mathfrak{H} )_*$ and $\rho$ is positive.  Then defining
$\vartheta _t$ and $\nu_t$ as in the proof of the last theorem we have
$\vartheta_t \geq \nu_t \geq 0$ for $t > 0$ and $\nu _t \rightarrow \hat
\pi (\rho )$ in norm as $t \rightarrow 0^+$.  Since $\vartheta_t - \nu_t
\geq 0$ and $\alpha _t(I) \leq I$ and $\pi$ is unital we have
\begin{align*}
\Vert\vartheta_t - \nu_t\Vert &= \vartheta_t(I) - \nu_t(I) = t^{-1}\Phi
(\rho )(\alpha_t(I) - I) + \Phi (\rho )(I) - \nu_t(I)
\\
&\leq \Phi (\rho )(I) - \nu_t(I) = \rho (I) - \nu_t(I) \rightarrow \rho (I)
- \hat \pi (\rho )(I) = 0.
\end{align*}
Hence, $\vartheta_t \rightarrow \hat \pi (\rho )$ in norm as $t \rightarrow
0^+$.  Since each $\rho \in \mathfrak{B} (\mathfrak{K} )_*$ is the linear
combination
of at most four positive elements of $\mathfrak{B} (\mathfrak{K} ) _*$ we have
$$
t^{-1}(\hat \alpha_t(\Phi (\rho )) - \Phi (\rho )) + \Phi (\rho )
\rightarrow \hat \pi (\rho )
$$
in norm as $t \rightarrow 0^+$.  Thus, $\Phi (\rho ) \in \mathfrak{D} (\hat
\delta )$ and $\hat \delta (\Phi (\rho )) + \Phi (\rho ) = \hat \pi (\rho
)$ for all $\rho \in \mathfrak{B} (\mathfrak{K} )_*$.  Hence, $\alpha$ is
derived from $\pi$.
\end{proof}

We will want to analyze the action of local cocycles on units.  Suppose
$\alpha$ is a unital $CP$-flow and $\alpha^d$ is the minimal dilation of
$\alpha$ to an $E_\subo $-semigroup acting on $\mathfrak{B} (\mathfrak{H}_1)$ as
described
in Theorem~\ref{thm:Bhat's-dilation-thm}.  A \empha{unit} for $\alpha^d$ is a
one parameter semigroup of isometries $V(t)$ which intertwine $\alpha^d$ so $(V(t)A =
\alpha_t^d(A)V(t)$ for all $A \in \mathfrak{B} (\mathfrak{H}_1)$ and $t \geq
0)$.
Units for $\alpha^d$ are in one to one correspondence with semigroups
$S(t)$ acting on $\mathfrak{H}$ with the property that the semigroup
$\Omega_t(A) = S(t)AS(t)^*$ is a trivially maximal subordinate of $\alpha$
(i.e., the mapping $A \rightarrow \alpha_t(A) - e^{st}\Omega_t(A)$ for $A
\in \mathfrak{B} (\mathfrak{H} )$ is completely positive for all $t \geq 0$
provided
$s \leq 0$ and the mapping is not positive for $s > 0$ and $t > 0)$.  The
next two theorems (see Theorems 4.46, 4.50 and 4.51 from \cite{BP})
describe such semigroups and the connection between them and units for the
dilated $E_\subo $-semigroup.

\begin{theorem}\label{thm:generator-pure-subordinate}%{Theorem 1.21}
Suppose $\alpha$ is a $CP$-flow over
$\mathfrak{K}$ and
$S(t)$ is a strongly continuous one parameter semigroup and $\Omega_t(A) =
S(t)AS(t)^*$ for $t \geq 0$ and $A \in \mathfrak{B} (\mathfrak{H} )$ is a
subordinate
of $\alpha .$  Then $S(t)$ is a strongly continuous one parameter semigroup
of contractions with generator $-D$ where $\mathfrak{D} (D) = \{f \in
\mathfrak{D}
(d^*): f(0) = Vf\}$ and $Df = -d^*f + cf$ where $c$ is a complex number
with non-negative real part and $V$ is a linear operator from $\mathfrak{H}$ to
$\mathfrak{K}$ with norm satisfying $\Vert V\Vert^2 \leq 2Re(c).  $ Furthermore,
if $\pi (A) = (2Re(c))^{-1}VAV^*$ for all $A \in \mathfrak{B} (\mathfrak{H} )$
and
$\gamma$ is the minimal $CP$-semigroup derived from $\pi$ then $\alpha$
dominates $\gamma .$  In the case $Re(c) = 0$ we define $\pi = 0.$

Conversely, if $c$ is a complex number with Re(c) $> 0$ and $V$ is a linear
operator from $\mathfrak{H}$ to $\mathfrak{K}$ with norm satisfying $\Vert
V\Vert ^2
\leq 2Re(c)$ and if $\pi (A) = (2Re(c))^{-1}VAV^*$ for $A \in \mathfrak{B}
(\mathfrak{H} )$ and $\gamma$ is the minimal $CP$-semigroup derived from $\pi$
and $\alpha$ dominates $\gamma$ then if $D$ is an operator with domain
$\mathfrak{D} (D) = \{f \in \mathfrak{D} (d^*): f(0) = Vf\}$ and $Df = -d^*f +
cf.$
Then $-D$ is the generator of a contraction semigroup $S(t)$ and if
$\Omega_t(A) =$ $S(t)AS(t)^*$ for $t \geq 0$ and $A \in \mathfrak{B}
(\mathfrak{H} )$
and $\alpha$ dominates $\Omega .$
\end{theorem}

\begin{theorem}\label{thm:trivially-max-subords-unit-dilatns}%{Theorem 1.22}
Suppose $\alpha$ is a unital $CP$-flow over
$\mathfrak{K}$ and $\alpha^d$ is the minimal dilation of $\alpha$ to an
$E_\subo $-semigroup
and suppose the relation between $\alpha$ and $\alpha^d$ is given by
$$
\alpha_t(A) = W^*\alpha_t^d(WAW^*)W
$$
for all $A \in \mathfrak{B} (\mathfrak{H} )$ (with $\mathfrak{H} = \mathfrak{K}
\otimes
L^2(0,\infty ))$ and $t \geq 0$ where $W$ is an isometry from $\mathfrak{H}$ to
$\mathfrak{H}_1$ and $WW^*$ is an increasing projection for $\alpha^d$ and
$\alpha^d$ is minimal over the range of $W.$ Then $\mathfrak{H} _1$ can be
expressed as $\mathfrak{H}_1 = \mathfrak{K}_1 \otimes L^2(0,\infty )$ and
$\alpha ^d$
is a $CP$-flow over $\mathfrak{K}_1$ so that if $U(t)$ and $U_1(t)$ are right
translation on $\mathfrak{H}$ and $\mathfrak{H} _1$ for $\alpha$ and $\alpha^d,$
respectively, then $U_1(t)W = WU(t)$ and $U_1(t)^*W = WU(t)^*$ for all $t
\geq 0.$   This means that $W$ as a mapping of $\mathfrak{H} = \mathfrak{K} \otimes
L^2(0,
\infty )$ into $\mathfrak{H}_1 = \mathfrak{K}_1 \otimes$ $L^2(0,\infty )$ can be
expressed in the form $W = W_1 \otimes I$ where $W_1$ is an isometry from
$\mathfrak{K}$ into $\mathfrak{K}_1.$

Suppose $S(t)$ is a strongly continuous semigroup of contractions of
$\mathfrak{H}$ and $\Omega$ given by $\Omega_t(A) = S(t)AS(t)^*$ for $A \in
\mathfrak{B}
(\mathfrak{H} )$ and $t \geq 0$ is a subordinate of $\alpha .$  Further assume
$\Omega$ is trivially maximal.  Then there is a unique strongly continuous
one parameter semigroup of isometries $S_1(t)$ which intertwine
$\alpha_t^d$ for each $t \geq 0$ and
$$
S(t) = W^*S_1(t)W
$$
for all $t \geq 0.$

Conversely, if $S_1(t)$ is a strongly continuous one parameter semigroup of
isometries which intertwine $\alpha_t^d$ for each $t \geq 0$ then if $S(t)$
is as defined in the equation above we have that $S(t)$ is a strongly continuous
one parameter semigroup of contractions so that $\Omega$ defined by
$\Omega_t(A) = S(t)AS(t)^*$ for $A \in \mathfrak{B} (\mathfrak{H} )$ and $t \geq
0$
is a subordinate of $\alpha$ which is trivially maximal.
\end{theorem}

We end this section with some notation and results which we will need in
the next section.  As we saw in Theorem~\ref{thm:bdry-wght-when-rep-normal}
the boundary weight map of the minimal $CP$-flow derived from $\pi$ is given by
$$
\omega (\rho ) = \hat \pi (\rho ) + \hat \pi (\hat \Lambda (\hat \pi (\rho
))) + \hat \pi (\hat \Lambda (\hat \pi (\hat \Lambda (\hat \pi (\rho )))))
+ \cdots  .
$$
We introduce some notation.  For $n \in \mathbb{N}$ we write $R_n(\phi )$ to
denote finite sum and $R(\phi )$ to denote the infinite series
$$
R_n(\phi ) = I + \phi + \phi^2 + \cdots + \phi^n \qquad \text{and} \qquad
R(\phi ) = I + \phi + \phi^2 + \cdots.
$$
Then the expression for $\omega$ above can be written
$$
\omega = \hat \pi R(\hat \Lambda\hat \pi ) = R(\hat \pi\hat \Lambda )\hat
\pi .
$$
Formally, $R(\phi ) = (I - \phi )^{-1}$, however the inverse in question
may not exist.  The sums above make sense in that the series
$$
\omega (\rho )(A) = \rho (\pi (A)) + \rho (\pi (\Lambda (\pi (A)))) +
\cdots
$$
converges absolutely for $A \in \mathfrak{A} (\mathfrak{H} )$.  This is seen by
setting $A = I - \Lambda$ and assuming $\rho \in \mathfrak{B} (\mathfrak{K} )_*$
is
positive.  As we saw in Theorem~\ref{thm:bdry-wght-when-rep-normal} the
series above for $\omega$ defines
the minimal $CP$-flow derived from $\pi .  $ We know from
Theorem~\ref{thm:restricted-weights} that the
truncated boundary weight map $\rho \rightarrow \omega\vert_t$ for $t > 0$
is the minimal $CP$-flow derived from the truncated boundary representation
$\phi_t^\#$.

%%%%%%%%%%%%%%%%%%% SECTION 2 %%%%%%%%%%%%%%%%%%%%%%%%%%%%%%%%%%%%%

\section{An Almost Type~I $CP$-Flow}\label{section:almost-typeI}

In this section we study $CP$-flows derived from a particular strongly
continuous $*$-representation $\pi$.  Let $\mathfrak{K}$ be the infinite tensor
product of $L^2(0,\infty )$ so $\mathfrak{K} = \otimes_{k=1}^\infty L^2(0,\infty
)$ with the reference vector (see \cite{vN} for details of infinite tensor
products of Hilbert spaces)
$$
F_\subo  = k_1 \otimes k_2 \otimes \cdots
$$
with
$$
k_i(x) = \lambda_ie^{-{\tfrac{1}{2}}\lambda_i^2x}
$$
for $x \geq 0$ where $\lambda_i > 0$ for $i = 1,2,\cdots$.  The Hilbert
space $\mathfrak{K}$ is spanned by product vectors of the form
$$
F = f_1 \otimes f_2 \otimes \cdots
$$
where
\begin{equation}\label{eq:inf-tensor-cond-1}
\sum_{i=1}^\infty \Vert f_i - k_i\Vert^2 < \infty %2.1
\end{equation}
The inner product between two such product vectors is given by
$$
(F,G) =\prod_{i=1}^\infty (f_i,g_i)
$$
We impose the following two conditions on the positive numbers $\lambda_i$.
\begin{equation}\label{eq:inf-tensor-cond-2}
\sum_{n=1}^\infty \lambda_n^{-2} < \infty\qquad \text{and} \qquad \sum_{n=1
}^\infty \frac {\vert\lambda_n-\lambda_{n+1}\vert^2}
{\lambda_n^2+\lambda_{n +1}^2} < \infty   %\tag2.2
\end{equation}
We note both these conditions are satisfied for $\lambda_n = n$ and the
second condition in not satisfied for $\lambda_n = 2^n$.  Let $W_\subo $ be the
unitary mapping of $\mathfrak{H} = \mathfrak{K} \otimes L^2(0,\infty )$ into
$\mathfrak{K}$ given by
\begin{equation}\label{eq:S-zero-formula}
S_\subo ((f_1 \otimes f_2 \otimes \cdots ) \otimes h) = h \otimes f_1 \otimes
f_2 \otimes \cdots % \tag2.3
\end{equation}
and let $\pi (A) = S_\subo AS_\subo ^*$ and $\Delta = e^{-x} \otimes e^{-x}
\otimes
\cdots $ where $e^{-x}$ is shorthand for the operation of multiplication by
$e^{-x}$ on $L^2(0,\infty )$.  The first sum condition insures that
$\Delta$ is not zero and the second condition insures that $S_\subo $ is well
defined.  Note $\pi$ is a normal $*$-representation of $\mathfrak{B}
(\mathfrak{H} )$
on $\mathfrak{B} (\mathfrak{K} )$.  Suppose $n$ is a positive integer.  We
define
$\mathfrak{K}_n$ as the tensor product of the Hilbert spaces $\mathfrak{B}
(L^2(0,\infty ))$ from $n+1$ on with the reference vector $F_{n\subo} =
\otimes_{i=k+n}^\infty k_i$.  Let $S_n$ be the linear mapping from
$\mathfrak{K}$ to $\mathfrak{K}_n$ which takes the product vector $F =
\otimes_{i=1}^ \infty
f_i \in \mathfrak{K}$ to the product vector $S_nF = \otimes_{i=1}^\infty f_i \in
\mathfrak{K}_n$.  From the second sum condition above one finds $S_n$ is well
defined and one checks that $S_n$ is unitary.  We define $Q_n(A) =
S_nAS_n^*$ for $A \in \mathfrak{B} (\mathfrak{K} )$.  Let
$\mathfrak{K}_n^{\prime} =
\otimes_{i=1}^nL^2(0,\infty )$.  We see that $\mathfrak{B} (\mathfrak{K} ) =
\mathfrak{B}
(\mathfrak{K}_n^{\prime}) \otimes \mathfrak{B} (\mathfrak{K} _n)$ and
$Q_m(I_n\otimes
Q_n(A)) = I_{n+m}\otimes Q_{n+m}(A)$ for $A \in \mathfrak{B} (\mathfrak{H} )$
where
$I_k$ is the unit in $\mathfrak{B} (\mathfrak{K}_k^{\prime})$ for $n,m,k =
1,2,\cdots$.

We have the formulae,
\begin{align*}
\pi (A \otimes A_\subo ) &= A_\subo  \otimes Q_1(A)
\\
\Lambda (A) &= A \otimes e^{-x}
\\
\pi (\Lambda (A)) &= e^{-x} \otimes Q_1(A)
\\
(\pi\Lambda )^n(A) &= e^{-x} \otimes e^{-x} \otimes \cdots  \otimes e^{-x}
\otimes Q_n(A)
\end{align*}
for $A \in \mathfrak{B} (\mathfrak{K} )$ and $n = 1,2,\cdots $ where there are
$n$
factors of $e^{-x}$ in the last equation.  For $A = A_1\otimes A_2 \otimes
\cdots $ we write these formulae
\begin{align*}
\pi ((A_1\otimes A_2\otimes \cdots ) \otimes A_\subo ) &= (A_\subo \otimes A_1
\otimes A_2 \otimes \cdots )
\\
\Lambda (A_1\otimes A_2\otimes \cdots ) &= (A_1 \otimes A_2 \otimes \cdots
) \otimes e^{-x}
\\
\pi (\Lambda (A_1\otimes A_2\otimes \cdots )) &= e^{-x} \otimes A_1 \otimes
A_2 \otimes \cdots
\\
(\pi\Lambda )^n(A_1\otimes A_1\otimes \cdots ) & = e^{-x} \otimes e^{-x}
\otimes \cdots \otimes e^{-x} \otimes A_1\otimes A_2\otimes \cdots
\end{align*}

We first note that $(\pi\Lambda )^n (I)$ converges to $\Delta$ as
$n \rightarrow \infty$.  We have
$$
(\pi\Lambda )^n(I) = e^{-x} \otimes e^{-x} \otimes \cdots  \otimes e^{-x}
\otimes I \otimes I \otimes \cdots
$$
where there are $n$ factors of $e^{-x}$ and we see that $(\pi\Lambda
)^n(I)$ forms a decreasing sequence of positive operators which must
converge strongly to a limit which is
$$
\Delta = e^{-x} \otimes e^{-x} \otimes \cdots .
$$
As we have mentioned the first sum condition on the $\lambda_n$ insures
that $\Delta$ is not zero.

Next we note that if $\pi (\Lambda (A)) = A$ then $A$ is a multiple of
$\Delta$ (i.e.  $A = c\Delta$ with $c \in \mathbb{C} )$.  Note that if $\pi
(\Lambda (A)) = A$ then if $A = A_1 + iA_2$ where $A_1$ and $A_2$ are
hermitian then $\pi (\Lambda (A_i)) = A_i$ for $i = 1,2$.  So it is enough
to show that if $A = A^*$ and $\pi (\Lambda (A)) = A$ then $A = \lambda
\Delta$ with $\lambda$ real.  Next note that if $A \in \mathfrak{B}
(\mathfrak{K} )$
is hermitian and $\pi (\Lambda (A)) = A$ and $\Vert A\Vert = 1$ then
$(\pi\Lambda )^n(I + A) \rightarrow \Delta + A$ as $n \rightarrow \infty$
and since $\Delta + A$ is the strong limit of positive operators we have
$\Delta + A$ is positive.  If $\Delta + A = \lambda\Delta$ it follows that
$A$ is a multiple of $\Delta $.  Hence, it is sufficient to show that if
$A$ is positive and $\pi (\Lambda (A)) = A$ then $A = \lambda\Delta$ with
$\lambda \geq 0$.  Suppose then that $A \in \mathfrak{B} (\mathfrak{K} )$ is
positive, $\Vert A\Vert = 1$ and $\pi (\Lambda (A)) = A$.  Since
$(\pi\Lambda )^n(I - A) \rightarrow \Delta - A \geq 0$ we have $0 \leq A
\leq \Delta$.  Recalling the reference vector $F_\subo $ we have
$$
(F_\subo ,\Delta F_\subo ) = (k_1,e^{-x}k_1)(k_2,e^{-x}k_2) \cdots   =
\frac{\lambda_1^2}{1+\lambda_1^2} \cdot \frac{\lambda_2^2}{1+\lambda_2^2} \cdots
$$
Since $\Delta \geq A \geq 0$ we have $(F_\subo ,AF_\subo ) = c(F_\subo ,\Delta
F_\subo )$ with
$c \in [0,1]$.  Now since $\pi (\Lambda (A)) = A$ if follows that
$$
A = e^{-x} \otimes Q_1(A) = e^{-x} \otimes e^{-x} \otimes Q_2(A) = \cdots
$$
and we have
$$
(\otimes_{i=n+1}^\infty k_i,Q_n(A) \otimes_{i=n+1}^\infty k_i) =
c\; \frac{\lambda _{n+1}^2}{1+\lambda_{n+1}^2} \cdot
\frac{\lambda_{n+2}^2}{1+\lambda_{n+2}^2} \cdots
$$
for $n = 1,2,\cdots$.  Now let
$$
F = \otimes_{i=1}^\infty f_i\qquad \text{and} \qquad G =
\otimes_{i=1}^\infty g_i
$$
be product vectors so that $f_i = g_i = k_i$ for $i \geq m$.  Then we see
that
\begin{align*}
(F,AG) &= (f_1,e^{-x}g_1)(f_2,e^{-x}g_2)\cdots (f_m,e^{-x}g_m)\; c\;
\frac{\lambda _{m+1}^2}{1+\lambda_{m+1}^2} \frac{\lambda_{m+2}^2}{1+\lambda_{m+2}^2}
\cdots\\
&= c(F,\Delta G).
\end{align*}
Since such vectors $F$ and $G$ are dense in $\mathfrak{K}$ we have $A = c\Delta
$.  Then we have proved the following lemma.

\begin{lemma}%{Lemma 2.1}
Suppose $\pi$ is the $*$-representation described
above.  Let $\Delta$ be as described above.  Then $\Delta =
\lim_{n\rightarrow\infty} ( \pi\Lambda )^n(I)$.  Furthermore, if $A \in
\mathfrak{B} (\mathfrak{K} )$ and $\pi (\Lambda (A)) = A$ then $A$ is a multiple
of
$\Delta$ (i.e., $A = c\Delta$ with $c \in \mathbb{C} )$.
\end{lemma}

We will need a stronger characterization of this property which is provided
by the following lemma.

\begin{lemma}\label{lemma:rho(delta)=0-trick}%{Lemma 2.2}
Suppose $\rho \in \mathfrak{B} (\mathfrak{K} )_*$ and
$\rho
(\Delta ) = 0$.  Then $\Vert (\hat \Lambda\hat \pi )^n(\rho )\Vert
\rightarrow 0$ as $n \rightarrow \infty$.
\end{lemma}
\begin{proof}  Suppose $\rho \in \mathfrak{B} (\mathfrak{K} )_*$ and $\rho
(\Delta ) =
0$.  Suppose $\epsilon > 0$.  Since $\rho$ can be approximated arbitrarily
well in norm by a finite sum of functionals $\rho_i$ of the form $\rho_i(A)
= (F_i,AG_i)$ with $F_i,G_i \in \mathfrak{K}$ for $i = 1,\cdots ,n$ and the
vectors $F_i$ and $G_i$ can be approximated by vectors $F_i^{\prime}$ and
$G_i^{\prime}$ which are finite sums of product vectors of the form $f_1
\otimes f_2 \otimes \cdots $ with $f_i = k_i$ for $i > m$ with $m$ some
large integer it follows that there is a functional $\eta$ so that
$\Vert\rho - \eta\Vert < {\tfrac{1}{2}}\epsilon (F_\subo ,\Delta F_\subo )$ and
$$
\eta (A) = \sum_{i=1}^n (F_i,AG_i)
$$
and each of the vectors $F_i$ and $G_i$ are of the form
$$
F \otimes (\otimes_{i=m+1}^\infty k_i)
$$
(i.e. they consist of sums of product vectors with factors $f_i = k_i$ for
$i > m)$.  Since we have
$$
\vert\eta (\Delta )\vert = \vert\rho (\Delta ) - \eta (\Delta )\vert \leq
\Vert \rho - \eta\Vert < {\tfrac{1}{2}}\epsilon (F_\subo ,\Delta F_\subo )
$$
Now let $\mu (A) = \eta (A) - (F_\subo ,AF_\subo )\eta (\Delta )(F_\subo ,\Delta
F_\subo )^{-1}.  $ Note
$$
\Vert\rho - \mu\Vert \leq \Vert\rho - \eta\Vert + \Vert\eta - \mu\Vert \leq
{\tfrac{1}{2}}\epsilon (F_\subo ,\Delta F_\subo ) + {\tfrac{1}{2}}\epsilon <
\epsilon .
$$
Since $(\hat \Lambda\hat \pi )^k(\mu )(A) = \mu ((\pi\Lambda )^k(A))$ and
$( \pi\Lambda )^k(A)$ is of the form
$$
(\pi\Lambda )^k(A) = e^{-x} \otimes e^{-x} \otimes \cdots  \otimes e^{-x}
\otimes A_m
$$
where there are $k$ factors of $e^{-x}$, and it follows from the form of $\mu$
that $\mu ((\pi\Lambda )^k(A)) = 0$ for $k \geq m$.  Hence, $\Vert (\hat
\Lambda \hat \pi )^k(\rho )\Vert < \epsilon$ for $k \geq m$ and we have
$\Vert (\hat \Lambda\hat \pi )^k(\rho )\Vert \rightarrow 0$ as $k
\rightarrow \infty$.
\end{proof}

Let $\alpha^1$ be the minimal $CP$-flow derived from $\pi$. If $\rho
\rightarrow \omega^1(\rho )$ is the boundary weight map for $\alpha^1$ then
$$
\omega^1(\rho ) = \hat \pi (\rho ) + \hat \pi\hat \Lambda\hat \pi (\rho ) +
\hat \pi\hat \Lambda\hat \pi\hat \Lambda\hat \pi (\rho ) + \cdots  =  \hat
\pi R(\hat \Lambda\hat \pi )
$$
where the short hand $R(\psi ) = I + \psi + \psi^2 + \cdots $ was
introduced in the last section.  We analyze $CP$-flows derived from $\pi$.
We begin with the following observation.

\begin{theorem}\label{thm:bdry-wght-for-derived-flow}%{Theorem 2.3}
Suppose $\alpha$ is a $CP$-flow derived from $\pi$
and $\omega$ is the boundary weight map for $\alpha$.  Then $\omega$ is of
the form
$$
\omega (\rho ) = \omega^1(\rho ) + \rho (\Delta )\xi
$$
for $\rho \in \mathfrak{B} (\mathfrak{K} )_*$ where $\omega^1$ is the boundary
weight map for $\alpha^1$ the minimal $CP$-flow derived from $\pi$ and $\xi \in
\mathfrak{A} (\mathfrak{H} )_*$ is a positive boundary weight on $\mathfrak{A}
(\mathfrak{H} )
= (I - \Lambda )^{\tfrac{1}{2}}\mathfrak{B} (\mathfrak{H} )(I - \Lambda
)^{\tfrac{1}{2}}$ with
$\xi (I - \Lambda ) \leq 1$ and $\alpha$ is unital (i.e. $\alpha_t(I) = I$
for $t \geq 0)$ if and only if $\xi (I - \Lambda ) = 1$.
\end{theorem}
\begin{proof}  Assume the hypothesis and notation of the theorem. Since
$\alpha$ is derived from $\pi$ we have $\omega (\rho - \hat \Lambda\hat \pi
(\rho )) = \hat \pi (\rho )$ for $\rho \in \mathfrak{B} (\mathfrak{K} )_*$.
Suppose
$\rho \in \mathfrak{B} (\mathfrak{K} )_*$.  Let $\rho_n = \rho + \hat
\Lambda\hat \pi
(\rho ) + \cdots  + (\hat \Lambda\hat \pi )^n(\rho )$.  Then we have
$$
\omega (\rho - (\hat \Lambda\hat \pi )^{n+1}(\rho )) = \hat \pi (\rho ) +
\hat \pi\hat \Lambda\hat \pi (\rho ) + \cdots  + \hat \pi (\hat \Lambda\hat
\pi )^n(\rho ).
$$
Now suppose $\rho (\Delta ) = 0$.  Then by Lemma~\ref{lemma:rho(delta)=0-trick}
we have $\Vert (\hat
\Lambda \hat \pi )^n\rho\Vert \rightarrow 0$ as $n \rightarrow \infty$ so
we have taking the limit as $n \rightarrow \infty$ that $\omega (\rho ) =
\omega^1(\rho )$ for $\rho \in \mathfrak{B} (\mathfrak{K} )_*$ with $\rho (\Delta
) =
0$.  Now suppose $\eta \in \mathfrak{B} (\mathfrak{K} )_*$ is positive and
$\eta(\Delta ) = 1$.  Then for arbitrary $\rho \in \mathfrak{B} (\mathfrak{K} )_*$
we have
$$
\omega (\rho ) = \omega (\rho - \rho (\Delta )\eta ) + \rho (\Delta )\omega
(\eta )  = \omega^1(\rho ) + \rho (\Delta )(\omega (\eta ) - \omega^1(\eta
)).
$$
Setting  $\xi = \omega (\eta ) - \omega^1(\eta )$ we have $\omega$ given in
terms of $\omega^1$ and $\xi$ as stated in the theorem.  Next we show $\xi$
is a positive.  Suppose $\rho \in \mathfrak{B} (\mathfrak{K} )_*$ is positive
and
$\rho (\Delta ) = 1$.  Then we have
$$
\omega ((\hat \Lambda\hat \pi )^n(\rho )) = \omega^1((\hat \Lambda\hat \pi
)^n)(\rho ) + \xi
$$
for each $n = 1,2,\cdots $ and since $\omega^1((\hat \Lambda\hat \pi
)^n)(\rho ) \rightarrow 0$ as a weight and since $(\hat \Lambda\hat \pi
)^n(\rho )$ is positive we have $\xi$ is the limit of positive weights so
$\xi$ is positive.  For $\rho \in \mathfrak{B} (\mathfrak{K} )_*$ we have
$$
\omega (\rho )(I - \Lambda ) = \omega^1(\rho )(I - \Lambda ) + \rho (\Delta
)\xi (I - \Lambda )
$$
and calculating $\omega^1(\rho )(I - \Lambda )$ we find
\begin{align*}
\omega^1(\rho )(I - \Lambda ) &= \rho ((I - \pi (\Lambda )) + (\pi (\Lambda
) - (\pi\Lambda )^2(\Lambda )) + \cdots  )
\\
&= \rho (I) - \rho (\Delta ).
\end{align*}
Hence, we have $\omega (\rho )(I - \Lambda )  = \rho (I) - \rho (\Delta )(1
- \xi (I - \Lambda ))$ for $\rho \in \mathfrak{B} (\mathfrak{K} )_*$.  Since we
have
then inequality $\omega (\rho )(I - \Lambda ) \leq \rho (I)$ for positive
$\rho \in \mathfrak{B} (\mathfrak{K} )_*$ we find $\xi (I - \Lambda ) \leq 1$
and
$\omega (\rho )(I - \Lambda ) = \rho (I)$ if and only if $\xi (I - \Lambda
) = 1$.
\end{proof}

We remark that it was shown in Theorem 4.62 of \cite{BP} if $\nu$ is a
positive element of $\mathfrak{B} (\mathfrak{H} )_*$ with $\nu (I) \leq 1$ and
$\xi$
is of the form
$$
\xi = (1 - \nu (\Lambda (\Delta ))^{-1}R(\hat \pi\hat \Lambda )\nu
$$
then $\omega$ of the form given in Theorem~\ref{thm:bdry-wght-for-derived-flow}
is the boundary weight map of a $CP$-flow $\alpha$ is unital if and only if $\nu
(I) = 1$.  In a subsequent paper we find necessary and sufficient conditions on
$\xi$ that $\omega$ as given in the statement of the above theorem is the
boundary weight map of a $CP$-flow over $\mathfrak{K}$.  If $\xi$ satisfies
these conditions we say $\xi$ is \empha{$q$-positive}.  In a subsequent paper we show
that the above formula for $\xi$ can be generalized to positive $\Lambda
(\Delta )$-weights with $\nu (I - \Lambda (\Delta )) \leq 1$.  We also show
that there are more general $\xi$.   For this paper we simply note that
there are plenty of $q$-positive $\xi$ which yield unital $CP$-semigroups
$\alpha$.

%%%%%%%%%%%%%%%%%%%%%%%%%% SECTION 3 %%%%%%%%%%%%%%%%%%%%%%%%%%%%%%%%%%%

\section{Local Flow Cocycles}

In this section we study the local flow cocycles associated with the
$CP$-flows constructed in the previous section.  Suppose $\alpha$ is a
unital $CP$-flow and $\alpha^d$ is the minimal dilation of $\alpha$ to an
$E_\subo $-semigroup.  As we saw in
Theorem~\ref{thm:trivially-max-subords-unit-dilatns} then
the $\alpha^d$ is also a
$CP$-flow over $\mathfrak{K}_1$ and the translation $U(t)$ on the Hilbert space
$\mathfrak{H}$ on which $\alpha$ lives dilate to the translations $U^1(t)$ on
the Hilbert space $\mathfrak{H}_1$ on which $\alpha^d$ lives.  Recall $t
\rightarrow C(t)$ is a local cocyle for $\alpha^d$ $C$ is a cocycle and
$C(t)$ commutes with $\alpha_t^d(\mathfrak{B} (\mathfrak{H}_1))$ for all $t \geq
0$.
The cocycle $C$ is a flow cocycle if $C(t)U^1(t) = U^1(t)$ for all $t \geq
0$.  Just as each local unitary cocycle corresponds to a hyper maximal
corner from $\alpha$ to $\alpha$, each local unitary flow cocycle for
$\alpha^d$ the dilation of a $CP$-flow over $\mathfrak{K}$ corresponds
to a hyper maximal flow corner $\gamma$ from $\alpha$ to $\alpha$. Here a
flow corner from $\alpha$ to $\alpha$ is a corner so that the matrix
$\Theta$ in Definition~\ref{def:corner} is a $CP$-flow over $\mathfrak{K} \oplus
\mathfrak{K} $.
Theorems~\ref{thm:contractive-cocycles=eigencorners} and
\ref{thm:matrix-contractive-cocycles=eigencorners} of Section~\ref{section:background} of
this paper are valid if one
replaces the word $CP$-semigroup with $CP$-flow and cocycle with flow
cocycle (see Theorem 4.54 of \cite{BP}).  One ambiguity that occurs in
speaking of flow corners is the following.  When one says $\gamma$ is a
maximal flow corner do we mean $\gamma$ is maximal as a flow corner or
simply maximal as a corner.  In Lemma 4.55 of \cite{BP} it was shown that
if $\alpha$ and $\beta$ are $CP$-semigroups and $\gamma$ is a flow corner
from $\alpha$ to $\beta$ then $\alpha$ and $\beta$ are $CP$-flows.  It then
follows that the two notions of maximality are the same.

We mention one technical problem.  Suppose $\alpha^d$ is the dilation of
the $CP$-flow $\alpha$ and $t \rightarrow C(t)$ is a contractive local
cocycle and $C(t)U^1(t) = \exp(-zt)U^1(t)$ for $t > 0$ where $z$ is a
complex number with positive real part.  Let $C^{\prime}(t) =
\exp(zt)C(t)$.  Then $C^{\prime}$ is a local flow cocycle, however, it is
not clear that it is contractive so there may not be a flow corner
associated with it.  Fortunately, Theorem 4.61 of \cite{BP} shows that
$C^{\prime}$ is contractive so there is a local flow corner associated with
it.  This means every contractive local cocycle $C$ with the property just
stated is of the form $C(t) =$ $\exp(-zt)C^{\prime}(t)$ for $t \geq 0$
where $C^{\prime}$ is a flow cocycle and $z$ is a complex number with non
negative real part.

Here we introduce some notation which we will use through out this
section.  As in the last section $\pi$ is the $*$-representation of
$\mathfrak{B} (\mathfrak{H} )$ on $\mathfrak{B} (\mathfrak{K} )$ constructed in
the last section. We
denote by $\xi$ a $q$-positive (usually unital) boundary weight and by
$\alpha = \alpha^\xi$ the $CP$-flow derived from $\pi$ associated with
$\xi$ as described in the last section.  The boundary weight map for
$\alpha$ is
$$
\omega (\rho ) = \omega^1(\rho ) + \rho (\Delta )\xi
$$
for $\rho \in \mathfrak{B} (\mathfrak{K} )_*$ where $\omega^1 = R(\hat \pi\hat
\Lambda )\hat \pi$.  Recall that $q$-positive means that $\omega$ given
above is the boundary weight of a $CP$-flow over $\mathfrak{K}$.  As we
mentioned in the last section the complete characterization of such $\xi$
will be given in a subsequent paper but for now we simply remark there are
many $q$-positive $\xi$ as given in the previous section.

If $z$ is a complex number with $\vert z\vert \leq 1$ we denote by
\begin{equation}\label{eq:omega-z}
\omega^z = zR(z\hat \pi\hat \Lambda )\hat \pi = z\hat \pi R(z\hat
\Lambda\hat \pi ) = z\hat \pi + z^2\hat \pi\hat \Lambda\hat \pi + z^3\hat
\pi\hat \Lambda \hat \pi\hat \Lambda\hat \pi + \cdots %\tag3.a
\end{equation}
where the sum converges as a boundary weight since the sum converges for $z
= 1$ where all the terms are positive.

Next we introduce a family of one parameter semigroups of isometries which
intertwine $\alpha$.  For $z$ any complex number we denote by $U_z$ the one
parameter semigroup of isometries of $\mathfrak{H} = \mathfrak{K} \otimes
L^2(0,\infty )$
$$
U_z(t) = \exp(-tD_z)\qquad \text{where} \qquad D_z = -d^* + {\tfrac{1}{2}}\vert
z
\vert^2I
$$
for $t \geq 0$ and $d$ is the operation of differentiation defined in the
last section and the domain $\mathfrak{D} (D_z) = \{ f \in \mathfrak{D} (d^*):
f(0) =
zS_\subo f\}$ where $S_\subo $ is the unitary operator mapping $\mathfrak{H}$
into
$\mathfrak{K}$ defined by \eqref{eq:S-zero-formula} in the last section as
$$
S_\subo ((f_1 \otimes f_2 \otimes \cdots  ) \otimes f_0) = f_0 \otimes f_1
\otimes f_2 \otimes \cdots %\tag2.3
$$
for $f_i \in L^2(0,\infty )$ and the $f_i$ satisfy condition
\eqref{eq:inf-tensor-cond-1} and
$S_\subo $ defines $\pi$ in that $\pi (A) = S_\subo AS_\subo ^*$ for $A \in
\mathfrak{B}
(\mathfrak{H} )$.  Note $U_0 = U$ the standard right translation and $D_0 = d$.

Suppose $w$ and $z$ are complex numbers.  We show the $U_z$ are a one
parameter family of isometries and the covariance $c(w,z)$ of $U_w$ with
$U_z$ is given by
\begin{equation}\label{eq:covariance-function}
U_w(t)^*U_z(t) = \exp(c(w,z)t)I = \exp({\tfrac{1}{2}} (2\overline {w}z - \vert
w\vert^2 - \vert z\vert^2)t) I%\tag3.1
\end{equation}
for $t \geq 0$.  For $f \in \mathfrak{D} (D_w)$ and $g \in \mathfrak{D} (D_z)$
we
have
\begin{align*}
\frac {d} {dt} (U_w(t)f,U_z(t)g) &= (d^*U_w(t)f,U_z(t)g) +
(U_w(t)f,d^*U_z(t)g)
\\
&\qquad - {\tfrac{1}{2}} (\vert w\vert^2 + \vert z\vert^2)(U_w(t)f,U_z(t)g)
\end{align*}
for $t \geq 0$.  Now we have
\begin{align*}
(d^*U_w(t)f,U_z(t)g) + (U_w(t)f,d^*U_z(t)g) &= ((U_w(t)f)(0),(U_z(t)g)(0))
\\
&= (wS_\subo U_w(t)f,zS_\subo U_z(t)g)
\\
&= \overline {w}z(U_w(t)f,U_z(t)g)
\end{align*}
for $t \geq 0$ where we have used the relation between $h(0)$ and $S_\subo h$
for $h$ in $\mathfrak{D} (D_z)$ or $\mathfrak{D} (D_w)$ and the fact that
$S_\subo $
is an
isometry.  Hence we have
$$
\frac {d} {dt} (U_w(t)f,U_z(t)g) = c(w,z)(U_w(t)f,U_z(t)g)
$$
for $t \geq 0$ and since the domains $\mathfrak{D} (D_w)$ and $\mathfrak{D}
(D_z)$
are dense in $\mathfrak{H}$ (see the argument in Lemma 4.44 of \cite{BP})
equation~\eqref{eq:covariance-function} follows.

Next we note that $S$ is a one parameter semigroup, so $\Omega$ given by
$$
\Omega_t(A) = S(t)^*AS(t)
$$
for $A \in \mathfrak{B} (\mathfrak{H} )$ and $t \geq 0$  satisfies $\alpha
\geq
\Omega$ (meaning the mapping A $\rightarrow \alpha_t(A) - \Omega_t(A)$ is
completely positive for $A \in \mathfrak{B} (\mathfrak{H} )$ and $t \geq 0)$ if
and
only if there are complex numbers $y,z$ with $Re(y) \geq 0$ so that
$$
S(t) = e^{-yt}U_z(t)
$$
for $t \geq 0$.  This follows from Theorem~\ref{thm:generator-pure-subordinate}
of Section~\ref{section:background} once one notes that the condition of the theorem is
satisfied if and only if the
mapping $A \rightarrow \pi (A) - (2Re(c))^{-1}VAV^*$ is completely positive
and since $\pi (A) = S_\subo AS_\subo ^*$ for $A \in \mathfrak{B} (\mathfrak{H}
)$ this is
the
case if and only if $V$ is an appropriate multiple of $S_\subo $.

Suppose $z \in \mathbb{C}$.  We show $U_z$ intertwines $\alpha$.  From the
result just established we have the mapping
$$
\beta_t(A) = \alpha_t(A) - U_z(t)AU_z(t)^*
$$
for $A \in \mathfrak{B} (\mathfrak{H} )$ and $t \geq 0$ is completely positive.
Suppose $t > 0$.  Then since $\alpha$ is unital we have we have $\beta_t(I)
= I - U_z(t)U_z(t)^*.  $ Since $U_z(t)$ is an isometry and $\beta_t$ is
positive we have
$$
0 \leq \beta_t(A) \leq I - U_z(t)U_z(t)^*
$$
for $A \in \mathfrak{B} (\mathfrak{H} )$ with $0 \leq A \leq I$ and consequently
$$
\beta_t(A) = (I - U_z(t)U_z(t)^*)\beta_t(A)(I - U_z(t)U_z(t)^*)
$$
and by linearity this extends to all $A \in \mathfrak{B} (\mathfrak{H} )$.  Then
we
have
$$
\alpha_t(A) = U_z(t)AU_z(t)^* + (I - U_z(t)U_z(t)^*)\alpha_t(A)(I - U_z(t)U
_z(t)^*)
$$
for all $A \in \mathfrak{B} (\mathfrak{H} )$.  And multiplying the above
equation on
the right by $U_z(t)$ we obtain $U_z(t)A = \alpha_t(A)U_z(t)$ so $U_z$
intertwines $\alpha $.  Summarizing our results to this point we have the
mapping $A \rightarrow \alpha_t(A) - V(t)AV(t)^*$ for $A \in \mathfrak{B}
(\mathfrak{H} )$ and $t \geq 0$ is completely positive where $V$ is a one
parameter
semigroup of contractions then the $V(t)$ are in fact multiples of a
semigroup $U_z$ of isometries which intertwine $\alpha$.

Now suppose $\alpha^d$ is the dilation of $\alpha$ to an $E_\subo $-semigroup on
$\mathfrak{H} _1$ as described in
Theorem~\ref{thm:trivially-max-subords-unit-dilatns}.  Then from
Theorem~\ref{thm:trivially-max-subords-unit-dilatns} we see
the mapping
\begin{equation}\label{eq:action-on-U1}
W^*U_z^1(t)W = U_z(t)%\tag3.2
\end{equation}
for $t \geq 0$ and $z \in \mathbb{C}$ give us a bijection from the $U_z$ to
intertwining semigroups of isometries $U_z^1$ which intertwine $\alpha^d$
where the covariance for the $U_z^1$ is the same as the covariance for the
$U_z$ given in equation~\eqref{eq:covariance-function}.  Also every
intertwining semigroup for $\alpha^d$ is of the form $V^1(t)$ $= e^{-yt}U_z^1(t)$
with $y,z \in \mathbb{C}$.

Next we describe the action of local cocycles on the units $U_z^1$.  One
checks that if $t \rightarrow C(t)$ is a local cocycle for $\alpha^d$ then
$C(t)U_z^1(t)$ is a intertwining semigroup for $\alpha^d$.  Now the action
of the local cocycles on the units $U_z^1$ must be the same as the action
of the local unitary cocycles on an $E_\subo $-semigroup of type~I$_1$.  The
local unitary cocycles generate automorphisms of the product systems
associated with an $E_\subo $-semigroup and these have been computed by Arveson
in section 3.8 of \cite{A3}.  Also, they have been computed by Bhat in
\cite{Bh}.  The general contractive flow local cocycles for a $CP$-flow of
type~I $CP$-semigroup are characterized in Theorem 2.11 in \cite{APP}.  Now
we characterize the action of the contractive local cocycles on units.  If
$C$ is a contractive local cocycle for $\alpha^d$ then there are complex
numbers $a,b,c,y \in \mathbb{C}$ with $\vert a\vert \leq 1$ and $Re(y) \geq 0$
so that the action of $C$ on the units $U_z^1$ is given by
\begin{equation}\label{eq:action-of-cocycle-on-units}
C(t)U_z^1(t) = \exp(t(-y - {\tfrac{1}{2}}\vert v + z\vert^2(1 - \vert a\vert^2) +
iIm(\overline {c}z))) U_{az+b}^1(t)%\tag3.3
\end{equation}
for $t \geq 0$ with
$$
v = - (1 - \vert a\vert^2)^{-1} (\overline {a}b + c)
$$
and when $\vert a\vert = 1$ then numbers $a,b,c \in \mathbb{C}$ above satisfy
the additional constraint $ac + b = 0$ so
$$
C(t)U_z^1(t) = e^{-t(y+iIm(a\overline{b}z))}U_{az+b}^1(t)
$$
The action of $C^*$ is obtained by making the replacements
$$
a \rightarrow \overline {a},\qquad b \leftrightarrow c\qquad \text{and }
\qquad y \rightarrow \overline {y}.
$$
In the case when $\vert a\vert = 1$ we parameterize $C$ with complex
numbers $(y,a,b)$ not using $c$ so the action of $C^*$ in this case is
given by
$$
C(t)^*U_z^1(t) = e^{-t(\overline{y}-iIm(\overline{b}z))}
U^{1}_{\overline{a}(z-b)}(t)
$$
If the cocycle is isometric then
$$
\vert a\vert = 1,\qquad ac + b = 0,\qquad \text{and}\quad  Re(y)= 0
$$
If the cocycle is a flow cocycle then $b = c = y = 0$ so the action of a
flow cocycle on the units $U_z^1$ is given by
$$
C(t)U_z^1(t) = e^{-{\tfrac{1}{2}}\vert z\vert^2(1-\vert a\vert^2)}U_{az}^1(t)
$$
for $t \geq 0$ and $z \in \mathbb{C}$.

If $C$ and $C^{\prime}$ are contractive local cocycles whose action on the
units is characterized by the $n$-tuples $(a,b,c,y)$ and $(a',b',c',y')$ as
describe above then the corresponding numbers for the product cocycle $t
\rightarrow C(t)C^{\prime}(t)$ are
$$(aa^{\prime},ab^{\prime} + b,\overline {a}^{\prime}c + c^{\prime},y +
y^{\prime} + iIm(\overline {c}b^{\prime}) - {\tfrac{1}{2}} r)
$$
where $r = 0$ if either $\vert a\vert = 1$ or $\vert a^{\prime}\vert = 1$
and otherwise
\begin{align*}
r &= (1-\vert a^{\prime}\vert^2)^{-1}\vert\overline
{a}^{\prime}b^{\prime}+c^{\prime} \vert^2 + (1-\vert a\vert^2)^{-1}\vert
b^{\prime}(1-\vert a\vert^2)-\overline {a}b-c\vert^2
\\
&\qquad - (1-\vert aa^{\prime}\vert^2)^{-1} \vert\overline
{aa}^{\prime}(ab^{\prime} +b) + \overline {a}^{\prime}c + c^{\prime}\vert^2
\end{align*}
is a nonnegative real function of
$(a,b,c,a^{\prime},b^{\prime},c^{\prime}).  $ Given the complexity of the
function $r$ above we wonder if there is a better parameterization of the
action of the local cocycles on the units.  If either of the local cocycles
above is unitary the number $r$ above is zero so the parameterization of
contractive local cocycles is much more difficult than the parameterization
of the unitary local cocycles.

We caution the reader that action of a local cocycle on the units $U_z^1$
does not completely determine the cocycle since in our case $\alpha^d$ is
not completely spatial.  In the next theorem we characterize the
contractive local flow cocycles which as we have explained is equivalent to
determining the flow corners from $\alpha$ to $\alpha $.  First we prove
the following lemma.

\begin{lemma}\label{lemma:bdry-reps-and-flow-corners}%{Lemma 3.1}
Suppose $\xi$ is a unital $q$-positive boundary
weight on $\mathfrak{A} (\mathfrak{H} )$ and $\alpha$ is the $CP$-flow over
$\mathfrak{K}$
derived from $\pi$ associated with $\xi$.  Suppose $\gamma$ is a flow
corner from $\alpha$ to $\alpha$ which means that
$$\Theta_t(
\left[\begin{matrix} A_{11}&A_{12}
\\
A_{21}&A_{22}
\end{matrix} \right]) =
\left[\begin{matrix} \alpha_t(A_{11})&\gamma_t(A_{12})
\\
\gamma_t^*(A_{21})&\alpha_t(A_{22})
\end{matrix} \right]$$
for $t > 0$ and $A_{ij} \in \mathfrak{B} (\mathfrak{H} )$ for $i,j = 1,2$ is a
$CP$-flow over $\mathfrak{K} \oplus \mathfrak{K}$.  Then there is a complex
number
$z$ with $\vert z\vert \leq 1$ so that $\Theta$ is derived from $\Pi_z$
given by
$$\Pi_z(
\left[\begin{matrix} A_{11}&A_{12}
\\
A_{21}&A_{22}
\end{matrix} \right]) =
\left[\begin{matrix} \pi (A_{11})&z\pi (A_{12})
\\
\overline{z}\pi (A_{21})&\pi (A_{22})
\end{matrix} \right] $$
for $A_{ij} \in \mathfrak{B} (\mathfrak{H} )$ for $i,j = 1,2$.  Furthermore, for
each
$w \in \mathbb{C}$ we have
$$
U_{zw}(t)A = e^{{\tfrac{1}{2}} t\vert w\vert^2(1-\vert
z\vert^2)}\gamma_t(A)U_w(t)
$$
and
$$
U_w(t)A = e^{{\tfrac{1}{2}} t\vert w\vert^2(1-\vert z\vert^2)}\gamma_t^*(A)U_w(t)
$$
for $A \in \mathfrak{B} (\mathfrak{H} )$ and $t \geq 0$.
\end{lemma}
\begin{proof}  Assume the hypothesis and notation of the theorem.  Let
$\alpha ^d$ and $\Theta^d$ be the dilation of $\alpha$ and $\Theta$ to
$E_\subo $-semigroups on $\mathfrak{H}_1$ and $\mathfrak{H}_1 \oplus
\mathfrak{H}_1$
and the
relation between the $CP$-flow and the dilated $E_\subo $-semigroup is as
described in Section~\ref{section:background} so
$$
\alpha_t(A) = W^*\alpha_t^d(WAW^*)W
$$
for $t \geq 0$ and $A \in \mathfrak{B} (\mathfrak{H} )$.  We will show that
there is
$z \in \mathbb{C}$ with $\vert z\vert \leq 1$ so that $\Theta$ is derived from
$\Pi_z$ as defined above.

First note that $U(t) \oplus U(t)$ intertwines $\Theta$.  Using this we
find the boundary representation of $\Theta$ is of the form
$$\Pi(A) = \Pi(
\left[\begin{matrix} A_{11}&A_{12}
\\
A_{21}&A_{22}
\end{matrix} \right]) =
\left[\begin{matrix} \pi (A_{11})&\phi (A_{12})
\\
\phi (A_{21})&\pi (A_{22})
\end{matrix} \right] $$
for $A$ in the domain of the generator of $\Theta$.   Since $\pi$ is pure
meaning the only subordinates of $\pi$ are of the form $\lambda\pi$ with $0
\leq \lambda \leq 1$ and $\Pi$ is completely positive it follows that $\phi
= z\pi$ for some $z \in \mathbb{C}$ with $\vert z\vert \leq 1$.  Note in
general the boundary representation is the direct sum of a normal and a
non normal representation of the domain of the generator but in our case we
are assured that the there is no non normal part because $\pi$ is unital
and therefore $\Pi$ is normal.  Thus the boundary representation of
$\Theta$ is $\Pi$ so $\Theta$ is derived from $\Pi $.

As we have seen since $\gamma$ is a flow corner from $\alpha$ to $\alpha$
there is a unique contractive local flow cocycle $C$ for $\alpha^d$ so that
$$
\gamma_t(A) = W^*C(t)\alpha_t^d(WAW^*)W
$$
for all $t \geq 0$ and $A \in \mathfrak{B} (\mathfrak{H} )$.  Then as we have
seen
there is number $y \in \mathbb{C}$ with $\vert y\vert \leq 1$ so that
$$
C(t)U_w^1(t) = \exp(-{\tfrac{1}{2}} t\vert w\vert^2(1 - \vert
y\vert^2))U_{yw}^1(t)
$$
for $t \geq 0$ and $w \in \mathbb{C}$.  Then we have

\begin{align*}
\gamma_t(A)U_w(t) &= W^*C(t)\alpha_t^d(WAW^*)WU_w(t) \\
&= W^*C(t)\alpha_t^d(WAW^*)U_w^1(t)W \\
&= W^*C(t)U_w^1WAW^*W \\
& = \exp(-{\tfrac{1}{2}} t\vert w\vert^2(1 - \vert
y\vert^2))W ^*U_{yw}^1(t)WA
\\
&= \exp(-{\tfrac{1}{2}} t\vert w\vert^2(1 - \vert y\vert^2))U_{yw}(t)A
\end{align*}
for $t \geq 0,\medspace w \in \mathbb{C}$ and $A \in \mathfrak{B} (\mathfrak{H}
)$.
Also
since $C$ is a local cocycle we have
$$
\gamma_t^*(A) = W^*\alpha_t^d(WAW^*)C(t)^*W = W^*C(t)^*\alpha_t^d(WAW^*)W
$$
so
$$
\gamma_t^*(A)U_w(t) = \exp(-{\tfrac{1}{2}} t\vert w\vert^2(1 - \vert
y\vert^2))U_w(t)A
$$
for $t \geq 0,\medspace w \in \mathbb{C}$ and $A \in \mathfrak{B} (\mathfrak{H}
)$.
Hence, we have proved the lemma provided we can show $y = z$.

We show $y = z$.  Let $d_2 = d \oplus d$ so $d_2$ is the ordinary
differential operator $d/dx$ on $\mathfrak{H} \oplus \mathfrak{H}$.  We use
capital
letters $F$ and $G$ to denote elements of $\mathfrak{H} \oplus \mathfrak{H}$ and
lower case letters $f,g$ to denote elements of $\mathfrak{H}$.  Recall from that
the boundary representation discussed in Section~\ref{section:background} for $\Theta$ is given by
$$
\Pi_z(A)F(0) = (AF)(0)
$$
for $F \in \mathfrak{D} (d_2^*)$ and $A \in \mathfrak{D} (\delta_2)$ where
$\delta_2$
is the generator of $\Theta$.  Suppose $w \in \mathbb{C}$ and $w \neq 0$.  Now
suppose $G = \{ 0,g\}$ and $g \in \mathfrak{D} (D_w)$ so $g \in \mathfrak{D}
(d^*)$
and $g(0) = wS_\subo g$.  Suppose $A \in \mathfrak{D} (\delta _2)$ and $A_{ij}
\in
\mathfrak{B} (\mathfrak{H} )$ are the matrix coefficients of $A$ for $i = 1,2$.
Now
from what we have shown we have
$$
\gamma_t(A_{12})U_w(t)g = \exp(-{\tfrac{1}{2}} t\vert w\vert^2(1 - \vert
y\vert^2))U _{yw}(t)A_{12}g
$$
for $t \geq 0$.  Since $-D_w$ is the generator of $U_w$ and $g \in \mathfrak{D}
(D_w)$ we have $U_w(t)g$ is differentiable in $t$ and since $A \in \mathfrak{D}
(\delta_2)$ we have $\gamma_t(A_{12})$ is differentiable in $t$ so the
expression on the left hand side of the above equation is differentiable in
$t$.  Hence, $U_{yw}(t)A_{12}g$ is differentiable in $t$ so $Ag \in \mathfrak{D}
(D_{yw})$ and we have $Ag \in \mathfrak{D} (d^*)$ and
\begin{align*}
(A_{12}g)(0) &= ywS_\subo A_{12}g = ywS_\subo A_{12}(w^{-1}S_\subo^*g(0))
\\
&= yS_\subo A_{12}S_\subo^*g(0) = y\pi (A_{12})g(0).
\end{align*}
Since $\Pi_z(A)F(0) = (AF)(0)$ we have
$$
(A_{12}g)(0) = z\pi (A_{12})g(0)
$$
and comparing the two equations we see $y = z$.
\end{proof}

\begin{theorem}\label{thm:flow-corners-and-bdr-wghts}%Theorem 3.2}
Suppose $\xi$ is a unital $q$-positive boundary
weight on $\mathfrak{A} (\mathfrak{H} )$ and $\alpha$ is the $CP$-flow over
$\mathfrak{K}$
derived from $\pi$ associated with $\xi$.  Suppose $\gamma$ is a flow
corner from $\alpha$ to $\alpha$ which means that
$$\Theta_t(
\left[\begin{matrix} A_{11}&A_{12}
\\
A_{21}&A_{22}
\end{matrix} \right]) =
\left[\begin{matrix} \alpha_t(A_{11})&\gamma_t(A_{12})
\\
\gamma_t^*(A_{21})&\alpha_t(A_{22})
\end{matrix} \right] $$
for $t > 0$ and $A_{ij} \in \mathfrak{B} (\mathfrak{H} )$ for $i,j = 1,2$ is a
$CP$-flow over $\mathfrak{K} \oplus \mathfrak{K}$ and if $\Omega$ is the
boundary
weight map for $\Theta$ then $\Omega$ is of the form
$$ \Omega(
\left[\begin{matrix} \rho_{11}&\rho_{12}
\\
\rho_{21}&\rho_{22}
\end{matrix} \right]) =
\left[\begin{matrix} \omega (\rho_{11})&\sigma (\rho_{12})
\\
\sigma^*(\rho_{21})&\omega (\rho_{22})
\end{matrix} \right] $$
for $\rho_{ij} \in \mathfrak{B} (\mathfrak{K} )_*$ for $i,j = 1,2$.  Then there
is a
unique complex number $z$ with $\vert z\vert \leq 1$ so if $z \neq 1$ then
$$
\sigma (\rho ) = \omega^z(\rho )
$$
for $\rho \in \mathfrak{B} (\mathfrak{K} )_*$ and if $z = 1$ then there is a
boundary
weight $\xi ^{\prime}$ so that
$$
\sigma (\rho ) = \omega^1(\rho ) + \rho (\Delta )\xi ^{\prime}
$$
for all $\rho \in \mathfrak{B} (\mathfrak{K} )_*$.
\end{theorem}
\begin{proof}  Assume the hypothesis and notation of the first paragraph of
the theorem.  Then from the previous lemma there is a unique $z \in \mathbb{C}$
with $\vert z\vert \leq 1$ so that $\Theta$ as given in the previous lemma
is derived from $\Pi_z$.  Since $\Theta$ is derived from $\Pi_z$ we have
repeating the argument of Theorem~\ref{thm:bdry-wght-for-derived-flow} that
$$
\sigma (\rho - z^{n+1}(\hat \Lambda\hat \pi )^{n+1}(\rho )) = z\hat \pi
(\rho ) + z^2\hat \pi\hat \Lambda\hat \pi (\rho ) + \cdots  + z^n\hat \pi
(\hat \Lambda \hat \pi )^n(\rho )
$$
Suppose $\rho (\Delta ) = 0$.  Then we have from
Lemma~\ref{lemma:rho(delta)=0-trick} that $\Vert
(\hat \Lambda\hat \pi )^n(\rho )\Vert \rightarrow 0$ as $n \rightarrow
\infty$ so we have
$$
\sigma (\rho ) = z\hat \pi R(z\hat \Lambda\hat \pi )(\rho )
$$
Choose a positive $\rho_1$ so that $\rho_1(\Delta ) = 1$ and we find
\begin{align*}
\sigma (\rho ) &= \sigma (\rho - \rho (\Delta )\rho_1) + \rho (\Delta
)\sigma (\rho_1)
\\
&= z\hat \pi R(z\hat \Lambda\hat \pi )(\rho ) + \rho (\Delta )(\sigma (\rho
_1) - z\hat \pi R(z\hat \Lambda\hat \pi )(\rho_1))
\end{align*}
Letting $\xi ^{\prime} = \sigma (\rho_1) - z\hat \pi R(z\hat \Lambda\hat
\pi )(\rho_1)$ we have
$$
\sigma (\rho ) = \omega^z(\rho ) + \rho (\Delta )\xi ^{\prime}
$$

Now since $\sigma$ is derived from $z\pi$ we have $\sigma (\rho - z\hat
\Lambda \hat \pi (\rho )) = z\hat \pi (\rho )$ for $\rho \in \mathfrak{B}
(\mathfrak{K} )_*$ and since $\omega^z$ is also derived from $z \pi$ we have the
same
equation is true for $\omega^z$ from which it follows that
$$
\rho (\Delta )\xi ^{\prime} - z\hat \Lambda\hat \pi\rho (\Delta )\xi
^{\prime} = (1 - z)\rho (\Delta )\xi ^{\prime} = 0
$$
for $\rho \in \mathfrak{B} (\mathfrak{K} )_*$.  For $z \neq 1$ the only solution
to
this equation is $\xi ^{\prime} = 0$.  Hence, if $z \neq 1$ we have
$$
\sigma (\rho ) = \omega^z(\rho ) = z\hat \pi (\rho ) + z^2\hat \pi\hat
\Lambda \hat \pi (\rho ) + z^3\hat \pi\hat \Lambda\hat \pi\hat \Lambda\hat
\pi (\rho ) + \cdots
$$
for $\rho \in \mathfrak{B} (\mathfrak{K} )_*$.
\end{proof}

\begin{theorem}%{Theorem 3.3}
Suppose $\xi$ is a unital $q$-positive boundary
weight on $\mathfrak{A} (\mathfrak{H} )$ and $\alpha$ is the $CP$-flow over
$\mathfrak{K}$
derived from $\pi$ associated with $\xi$.  Suppose $\alpha^d$ is the
minimal dilation of $\alpha$ to an $E_0$-semigroup of $\mathfrak{B}
(\mathfrak{H}_1)$
as given in Theorem~\ref{thm:cocycle-conj-dil=hypermax-corner} so
$$
\alpha_t(A) = W^*\alpha_t^d(WAW^*)W
$$
for $A \in \mathfrak{B} (\mathfrak{H} )$ and $t \geq 0$.  Then there is a
bijection
from the units $U_z^1$ of $\alpha^d$ onto the units of $U_z$ of $\alpha$
given by
$$
W^*U_z^1(t)W = U_z(t)
$$
for $t \geq 0$ and $z \in \mathbb{C}$.  Suppose $t \rightarrow C(t)$ is a local
unitary local cocycle which fixes $U_0^1$ so $C(t)U_0^1(t) = U_0^1(t)$ for
$t \geq 0$.  Then $C(t)U_z^1(t) = U_z^1(t)$ for $t \geq 0$ and all $z \in
\mathbb{C}$.  This means that the action of the local unitary cocycles on the
units contains no rotations.
\end{theorem}
\begin{proof}  Assume the hypothesis and notation of the theorem.  Let
$$
\gamma_t(A) = W^*C(t)\alpha_t^d(WAW^*)W
$$
for $A \in \mathfrak{B} (\mathfrak{H} )$ and $t \geq 0$.  From
Theorem~\ref{thm:cocycle-conj-dil=hypermax-corner} we have
$\gamma$ is a hypermaximal flow corner from $\alpha$ to $\alpha$.  Since
$C(t)$ is a unitary local cocycle which fixes $U_0^1$ we know from the
general properties of such cocycles discussed before
Lemma~\ref{lemma:bdry-reps-and-flow-corners} there is a
complex number $y$ of modulus one so that $C(t)U_w^1(t) = U_{yw}^1(t)$ for
all $t \geq 0$ and $w \in \mathbb{C} $.  Since $\gamma$ is a flow corner from
$\alpha$ to $\alpha$ we know that there is a complex number $z$ so that
$$
U_{zw}(t)A = e^{{\tfrac{1}{2}}\vert w\vert^2(1 - \vert z\vert^2)}
\gamma_t(A)U _w(t)
$$
for $w \in \mathbb{C} ,\medspace t \geq 0$ and $A \in \mathfrak{B} (\mathfrak{H}
)$.
Then we have
\begin{align*}
\gamma_t(A)U_w(t) &= W^*C(t)\alpha_t^d(WAW^*)WU_w(t) \\
& = W^*C(t)\alpha_t^d(WAW^*)U_w^1(t)W \\
&= W^*C(t)U_w^1(t)WAW^*W \\
& = W^*U_{yw}^1(t)WA = U_{yw}(t)A
\end{align*}
for all $w \in \mathbb{C} ,\medspace t \geq 0$ and $A \in \mathfrak{B}
(\mathfrak{H}
)$.
Comparing the two equations we see $y = z$ so $\vert z\vert = 1$.

To complete the proof to the theorem all we need do is show $z = 1$.
Suppose $\vert z\vert = 1$ and $z \neq 1$.  Now we apply
Theorem~\ref{thm:flow-corners-and-bdr-wghts}.  Let
$\Theta$ be the $CP$-flow described in the theorem.   Assuming the notation
of Theorem~\ref{thm:flow-corners-and-bdr-wghts} we have $\sigma (\rho ) =
\omega^z(\rho )$ for $\rho \in \mathfrak{B} (\mathfrak{K} )_*$.  We show this
implies that $C$ is not unitary.
From
Theorem~\ref{thm:cocycle-conj-dil=hypermax-corner} we know $C$ is a unitary
cocycle if and only if $\gamma$ is
hypermaximal.  But $\gamma$ is not hypermaximal as can be seen as follows.
Let $\Theta^1$ be the $CP$-semigroup of $\mathfrak{B} (\mathfrak{H} \oplus \mathfrak{H} )$ given by
$$\Theta_t^1(
\left[\begin{matrix} A_{11}&A_{12}
\\
A_{21}&A_{22}
\end{matrix} \right]) =
\left[\begin{matrix} \eta_t(A_{11})&\gamma_t(A_{12})
\\
\gamma_t^*(A_{21})&\eta_t(A_{22})
\end{matrix} \right] $$
for $t > 0$ and $A_{ij} \in \mathfrak{B} (\mathfrak{H} )$ for $i,j = 1,2$ where
$\eta$ is the minimal CP-flow derived from $\pi$.  Note the
boundary weight map weight map $\Omega^1$ for $\Theta^1$ is of the form
$$\Omega^1(
\left[\begin{matrix} \rho_{11}&\rho_{12}
\\
\rho_{21}&\rho_{22}
\end{matrix} \right]) =
\left[\begin{matrix} \omega^1(\rho_{11})&\omega^z(\rho_{12})
\\
\omega^{\overline{z}}(\rho_{21})&\omega^1(\rho_{22})
\end{matrix} \right] $$
for $\rho_{ij} \in \mathfrak{B} (\mathfrak{K} )_*$ for $i,j = 1,2$.  Now we see
$\gamma$ is not hypermaximal since $\alpha \geq \eta$ and if $\gamma$ were
hypermaximal we would have $\alpha = \eta$.  Hence, if $z \neq 1$ we have
$C$ is not unitary so the action of the local unitary cocycles does not
contain the rotations.
\end{proof}

%%%%%%%%%%%%%%%%%%%%%%%%% SECTION 4 %%%%%%%%%%%%%%%%%%%%%%%%%%%%%%%%%%%%

\section{Conclusion}

Here we present our conclusions.  Suppose $\mathfrak{K}$ is a separable Hilbert
space and $\mathfrak{H} = \mathfrak{K} \otimes L^2(0,\infty )$ and $S$ is a
unitary
mapping from $\mathfrak{H}$ onto $\mathfrak{K}$ and $\pi (A) = SAS^*$ for $A \in
\mathfrak{B} (\mathfrak{H} )$.  Note $\pi$ is an irreducible $*$-representation
of
$\mathfrak{B} (\mathfrak{H} )$ on $\mathfrak{B} (\mathfrak{K} )$.  Suppose
$\Lambda$ is the
mapping of $\mathfrak{B} (\mathfrak{K} )$ into $\mathfrak{B} (\mathfrak{H} )$
given by
$$
(\Lambda (A)F)(x) = e^{-x}AF(x)
$$
for $A \in \mathfrak{B} (\mathfrak{K} )$ and $x \geq 0$ for all
$\mathfrak{K}$-valued
function $F \in \mathfrak{H}$.  Let
$$
\Delta  = \lim_{n\rightarrow\infty} (\pi\Lambda )^n(I)
$$
where the limit exists in the sense strong convergence since the terms are
decreasing.  Assume $\Delta \neq 0$.  Assume further that for all $\rho \in
\mathfrak{B} (\mathfrak{K} )_*$ with $\rho (\Delta ) = 0$ we have $\Vert (\hat
\Lambda\hat \pi )^n(\rho )\Vert \rightarrow 0$ as $n \rightarrow \infty$.
It follows from this that if $\pi (\Lambda (A)) = A$ then $A =
\lambda\Delta$.

Then there are unital $CP$-semigroups $\alpha$ of $\mathfrak{B} (\mathfrak{H} )$
which are intertwined by the shifts $U$ and there boundary weight maps are
given by
$$
\omega (\rho ) = \omega^1(\rho ) + \rho (\Delta )\xi
$$
for $\rho \in \mathfrak{B} (\mathfrak{K} )_*$ where $\omega^1$ is the boundary
weight
map for the minimal $CP$-flow derived from $\pi$ and it is given by
$$
\omega^1(\rho ) = \hat \pi (\rho ) + \hat \pi\hat \Lambda\hat \pi (\rho ) +
\hat \pi\hat \Lambda\hat \pi\hat \Lambda\hat \pi (\rho ) + \cdots  =  \hat
\pi R(\hat \Lambda\hat \pi )
$$
and $\xi \in \mathfrak{A} (\mathfrak{H} )_*$ is a positive boundary weight on
$$
\mathfrak{A} (\mathfrak{H} ) = (I - \Lambda )^{\tfrac{1}{2}}\mathfrak{B}
(\mathfrak{H} )(I - \Lambda
)^{\tfrac{1}{2}}
$$
with $\xi (I - \Lambda ) \leq 1$ and $\alpha$ is unital (i.e. $\alpha_t(I)
= I$ for $t \geq 0)$ if and only if $\xi (I - \Lambda ) = 1$.  The boundary
weight $\xi$ satisfies certain positivity conditions which we analyze in a
separate paper.  It was shown in Theorem 4.62 of \cite{BP} that if $\nu$ is a
positive element of $\mathfrak{B} (\mathfrak{H} )_*$ with $\nu (I) \leq 1$ and
$\xi$
is of the form
$$
\xi = (1 - \nu (\Lambda (\Delta ))^{-1}R(\hat \pi\hat \Lambda )\nu
$$
then $\omega$ as given above is the boundary weight map of a $CP$-flow
$\alpha$ is unital if and only if $\nu (I) = 1$.

Then if $\alpha$ is such a unital $CP$-flow then $\alpha$ has a Bhat
dilation to an $E_0$-semigroup $\alpha^d$.  This $E_0$-semigroup is of
index one.  The action of the local unitary cocycles on the units for this
$E_\subo$-semigroup in not two-fold transitive.  The Hilbert space for the
dilation is of the form $\mathfrak{H}_1 = \mathfrak{K}_1 \otimes L^2(0,\infty )$
and
if $U^1(t)$ is right translation by $t$ on $\mathfrak{H}_1$ then $U^1$ is a unit
for $\alpha^d$ meaning
$$
U^1(t)A = \alpha_t^d(A)U^1(t)
$$
for all $A \in \mathfrak{B} (\mathfrak{H}_1)$ and $t \geq 0$.  If $C(t)$ is a
unitary
local cocycle for $\alpha^d$ so $C(t) \in \alpha_t^d(\mathfrak{B}
(\mathfrak{H}_1))$
for all $t \geq 0$ and the $C(t)$ are unitary operators satisfying the
relation $C(t)\alpha_t^d(C(s)) = C(t + s)$ for $s,t \geq 0$ and if
$C(t)U^1(t) = U^1(t)$ for $t \geq 0$ then $C(t) = I$ for all $t \geq 0$.
This means the action of the gauge group on the units of $\alpha^d$ is a
smaller group than then for an $E_0$-semigroup type~I$_1$.  Also, this means
$\alpha^d$ is not cocycle conjugate to the tensor product of a semigroup of
type~II$_0$ with an type~I$_1$ for if this was the case the action of gauge
group of the units would contain all the Euclidean transformations just as
the action for an $E_0$-semigroup of type~I$_1$. The same reasoning applies
to the $E_0$-semigroups of type II$_1$ corresponding to the examples of product
systems constructed by Tsirelson in \cite{Tsirelson2004}.

The full action for the gauge group on a type~I$_1$ is the Euclidean group
who action on $\mathbb{C}$ is given by $z \rightarrow az + b$ for $a,b \in
\mathbb{C}$ and $\vert a\vert = 1$.  In our examples we have the further
restriction $a = 1$.  Tsirelson has examples where there are the
restrictions $a = 1$ and Im(b) $= 0$.  It is quite possible in our case
there may be further restrictions.  It may be that $b$ lies on a one
dimensional line or even the further restriction $b = 0$.  This would be
interesting since it would give an example of an action which is rigid.
That means that if $C(t)$ is a local unitary cocycle and $U$ is an unit
then
$$
C(t)U(t) = e^{i\lambda t}U(t)
$$
for $t \geq 0$.

We are somewhat embarrassed to report that in order to establish this
result all that is required is to determine whether certain fairly simple
first order differential equations with constant coefficients have a
bounded solution or not.  The equation are parameterize by the complex
numbers $(a,b)$ with $\vert a\vert = 1$.   We have shown that if $a \neq 1$
the equations have no solution.  If the equations never have solutions the
action is rigid.  If the equations have solutions when $b$ lies on a one
dimensional line we are in the situation Tsirelson found and if the
equations have a solution for all $b$ then we are in the case where we have
transitivity of the gauge group on the units but no two fold transitivity.

As the reader can probably guess the feature that make these equation
interesting and difficult is that they involve infinitely many variables.
We will present them in a longer and more detailed paper.

%%%%%%%%%%%%%%%%%%%%%%%%%%%%%%%%%%%%%%%%%%%%%%%%
\providecommand{\bysame}{\leavevmode\hbox to3em{\hrulefill}\thinspace}
\bibliographystyle{amsplain}

\end{document}